\documentclass[12pt]{article}

\usepackage{amsmath}
\usepackage{amsthm}
\usepackage{amssymb}
\usepackage{amscd}

\theoremstyle{definition}
\newtheorem{theorem}{Theorem}[section]
\newtheorem{prop}[theorem]{Proposition}
\newtheorem{lemma}[theorem]{Lemma}
\newtheorem{corollary}[theorem]{Corollary}

\newtheorem{example}[theorem]{Example}
\newtheorem{remark}[theorem]{Remark}

\numberwithin{equation}{section}
\newenvironment{demo}[1]{%
 \trivlist
 \item[\hskip\labelsep
       {\bf #1.}]
}{%
\hfill\qedsymbol
 \endtrivlist
}

\newcommand\defterm[1]{{\sl #1}\/}
\newcommand\newterm[1]{{\sl #1}\/}

\newcommand\Int{\mathbb{Z}}
\newcommand\Rat{\mathbb{Q}}
\newcommand\ep{\varepsilon}
\newcommand\sgn{\operatorname{sgn}}
\newcommand\trans{{}^t\!}
\newcommand{\vectx}{X}

\newcommand\GL{\mathbf{GL}}    
\newcommand\Symp{\mathbf{Sp}}  
\newcommand\Orth{\mathbf{O}}   

\title{
A Compound Determinant Identity \\
for Rectangular Matrices
}

\author{
Masao ISHIKAWA\footnote{
Department of Mathematics,
Faculty of Education, University of the Ryukyus, Nishihara, Okinawa 903-0213, Japan,
{\tt ishikawa@edu.u-ryukyu.ac.jp}
}, \ 
Masahiko ITO\footnote{
School of Science and Technology for Future Life,
Tokyo Denki University, Tokyo 101-8457, Japan, {\tt mito@cck.dendai.ac.jp
}
} \ 
and 
Soichi OKADA\footnote{
Graduate School of Mathematics, Nagoya University, Furo-cho, Chikusa-ku, Nagoya 464-8602, Japan,
{\tt okada@math.nagoya-u.ac.jp}
}
}

\date{
}

\begin{document}

\maketitle

\begin{abstract}
A compound determinant identity for minors of rectangular matrices 
is established.
As an application, we derive Vandermonde type determinant formulae for 
classical group characters. 
\end{abstract}

\section{
Introduction
}

Cauchy (1812) \cite{C} and Sylvester (1851) \cite{S} established 
the following identity for an arbitrary square matrix 
(see \cite[pp.\,99--131 of vol.\,I, pp.\,193--197 of vol.\,II]{Mu} 
or \cite[pp.\,87--89]{T}):

\begin{prop} 
\label{prop:Sylv}
Let $s$ and $n$ be positive integers such that $s \geq n$.
Let $A=(a_{ij})_{1 \leq i, j \leq s}$ be an arbitrary square matrix 
of size $s$. 
We denote by $A^{I}_{J}$ the submatrix 
$(a_{i_{k},j_{l}})_{1 \leq k, l \leq n}$ of size $n$ 
for two subsets $I= \{ i_1 < \cdots < i_n \}$ and 
$J = \{ j_1 < \cdots < j_n \}$ of $[s] = \{ 1, 2, \dots, s \}$.
Then we have
\begin{equation}
\det \left( \det A^{I}_{J} \right)_{I,\,J}
=
( \det A )^{\binom{s-1}{n-1}},
\label{eq:Sylv}
\end{equation}
where the rows and columns of the determinant in the left-hand 
side are indexed by $n$-element subsets $I$, $J$ of $[s]$ and 
arranged increasingly in the lexicographic ordering.
\end{prop}

Using this proposition, in \cite{IO} the second and third authors 
gave a simple proof to the fact that 
a determinant formed by multiple ${}_{2r}\psi_{2r}$ 
basic hypergeometric series is evaluated as a product of 
$q$-gamma functions \cite{AI1, I}.
As limiting cases, the hypergeometric determinant formula 
includes the following determinant formulae for Schur functions
 $\GL_n(\lambda; x)$ 
and symplectic Schur functions $\Symp_{2n}(\lambda; x)$.
(Similar formulae hold for classical group characters. 
See \cite{IK}.) 

\begin{prop}
\label{prop:detS}
Let $s$ and $n$ be positive integers such that $s \geq n$ and 
let $\vectx = (x_1, \dots, x_s)$ be an $s$-tuple of indeterminates.
For an $n$-element subset $I = \{ i_1 < \cdots < i_n \}$, we put 
$\vectx_{I}=(x_{i_1}, \dots, x_{i_n})$.
Then we have
\begin{gather}
\det \Big(
\GL_n(\lambda;\vectx_I)
\Big)_{\lambda, \, I}
=
\Bigg(
\prod_{1 \le i < j \le s} (x_i - x_j)
\Bigg)^{\!\binom{s-2}{n-1}},
\label{eq:detS}
\\
\det \Big(
\Symp_{2n}(\lambda;\vectx_I)
\Big)_{\lambda, \, I}
=
\Bigg(
\prod_{1 \le j < k \le s}
 \frac{(x_k-x_j)(1-x_j x_k)}{x_jx_k}
\Bigg)^{\!\binom{s-2}{n-1}},
\label{eq:detchi}
\end{gather}
where the row index $\lambda$ runs over all partitions whose 
Young diagrams are contained in the rectangular diagram $((s-n)^n)$, 
and the column index $I$ runs over all $n$-element subset of $[s]$.
\end{prop}

The aim of this paper is to establish another formula 
(Theorem~\ref{th:main}) for compound determinant 
like the Cauchy--Sylvester identity. 
As a consequence, we can derive identities involving classical group 
characters (Theorem~\ref{th:Schur}) similar to Proposition~\ref{prop:detS}. 
The identities for Schur functions and symplectic Schur functions 
generalize the formulae given in \cite[Corollary~3.1 and Proposition~7.1]{AI2}, 
which are obtained as limiting cases of another hypergeometric 
determinant discussed in \cite{AI2}.

In order to report the main result, we introduce some notations 
and terminologies. 
For any integers $r$ and $s$, we use the symbols 
$[r,s] = \{ r, r+1, \dots, s \}$ and $[r] = [1,r]$.
If $r > s$, then we use the convention that $[r,s] = \emptyset$.
If $S$ is a finite set and $r$ is a nonnegative integer, 
let $\binom{S}{r}$ denote the set of all $r$-element subsets of $S$.
For $I = \{ i_1 < \cdots < i_r \}$ and $J = \{ j_1 < \cdots < j_r \}$ 
in $\binom{[N]}{r}$, we write $I < J$ if there is an index 
$k$ such that
$$
i_1 = j_1, \quad \dots, \quad  i_{k-1} = j_{k-1}, \quad  i_k < j_k.
$$
This gives a total order on $\binom{[N]}{r}$,
which is called the \defterm{lexicographic order}.

Let $s$ and $n$ be positive integers. 
Let $Z_{s,n}$ denote the set of (weak) compositions of $n$ 
with at most $s$ parts, i.e.,
$$
Z_{s,n}
=
\{
\mu=(\mu_{1},\dots,\mu_{s}) \in \Int^s
\,|\,
\mu_{1} \geq 0, \dots, \mu_{s} \geq 0, \,
\mu_{1} + \cdots + \mu_{s} = n
\},
$$
and let $Z^0_{s,n}$ denote the set of compositions of $n$ 
with exactly $s$ parts, i.e.,
$$
Z^0_{s,n}
=
\{
\mu=(\mu_{1},\dots,\mu_{s}) \in \Int^s
\,|\,
\mu_{1} > 0, \dots, \mu_{s} > 0, \,
\mu_{1} + \cdots + \mu_{s} = n
\}.
$$
Then the cardinalities of these sets are given by
$$
\# Z_{s,n} = \binom{s+n-1}{n},
\quad
\# Z^0_{s,n}
= \# Z^0_{s,n-s}
= \binom{s+(n-s)-1}{n-s}
= \binom{n-1}{n-s}.
$$
We introduce the following total order on  $Z_{s,n}$.
For $\lambda$, $\mu \in Z_{s,n}$, we write $\lambda < \mu$ 
if there is an index $k$ such that
$$
\lambda_1 = \mu_1, \quad \dots, \quad \lambda_{k-1} = \mu_{k-1},
\quad \lambda_{k}>\mu_{k}.
$$
To each integer sequence $\mu = (\mu_1, ..., \mu_s)$ with
$0 \le \mu_i \le n$ for all $i$,
we associate an $n$-element subset 
$\iota_{s,n}(\mu) \in \binom{[sn]}{n}$ defined by
\begin{equation}
\iota_{s,n}(\mu)
=
\bigsqcup_{i=1}^{s}
\left[ (i-1)n+1, (i-1)n+\mu_{i} \right].
\label{eq:iota}
\end{equation}
If $\iota_{s,n}$ is restricted to $Z_{s,n}$,
then 
it gives an injection $\iota_{s,n} \,:\, Z_{s,n} \rightarrow \binom{[sn]}{n}$, 
and one readily sees that, for $\lambda$, $\mu \in Z_{s,n}$,
$\iota_{s,n}(\lambda) < \iota_{s,n}(\mu)$ if and only if $\lambda < \mu$.
We may write $\iota=\iota_{s,n}$ for brevity
if there is no ambiguity on $s$ and $n$.

For example, if $s=3$ and $n=2$, then
$$
Z_{3,2}
= \{ (2,0,0), (1,1,0), (1,0,1), (0,2,0), (0,1,1), (0,0,2) \},
$$
and the map $\iota_{3,2}$ defined in (\ref{eq:iota}) is given by
\begin{gather*}
\iota_{3,2}(2,0,0) = \{ 1, 2 \}, \quad
\iota_{3,2}(1,1,0) = \{ 1, 3 \}, \quad
\iota_{3,2}(1,0,1) = \{ 1, 5 \},
\\
\iota_{3,2}(0,2,0) = \{ 3, 4 \}, \quad
\iota_{3,2}(0,1,1) = \{ 3, 5 \}, \quad
\iota_{3,2}(0,0,2) = \{ 5, 6 \}.
\end{gather*}
This injection $\iota_{3,2}$ can be visualized by the following picture:
$$
\begin{array}{ccccccc}
& &
\raisebox{-3pt}{
\setlength{\unitlength}{1.2pt}
\begin{picture}(60,10)
\put(0,0){\makebox(10,10){$1$}}
\put(10,0){\makebox(10,10){$2$}}
\put(20,0){\makebox(10,10){$3$}}
\put(30,0){\makebox(10,10){$4$}}
\put(40,0){\makebox(10,10){$5$}}
\put(50,0){\makebox(10,10){$6$}}
\end{picture}
}
&\qquad&
& &
\raisebox{-3pt}{
\setlength{\unitlength}{1.2pt}
\begin{picture}(60,10)
\put(0,0){\makebox(10,10){$1$}}
\put(10,0){\makebox(10,10){$2$}}
\put(20,0){\makebox(10,10){$3$}}
\put(30,0){\makebox(10,10){$4$}}
\put(40,0){\makebox(10,10){$5$}}
\put(50,0){\makebox(10,10){$6$}}
\end{picture}
}
\\
(2,0,0) & \mapsto &
\raisebox{-3pt}{
\setlength{\unitlength}{1.2pt}
\begin{picture}(60,10)
\put(0,10){\line(1,0){60}}
\put(0,0){\line(1,0){60}}
\put(0,0){\line(0,1){10}}
\put(10,0){\line(0,1){10}}
\put(20,0){\line(0,1){10}}
\put(30,0){\line(0,1){10}}
\put(40,0){\line(0,1){10}}
\put(50,0){\line(0,1){10}}
\put(60,0){\line(0,1){10}}
\put(0,0){\makebox(10,10){$\bigcirc$}}
\put(10,0){\makebox(10,10){$\bigcirc$}}
\end{picture}
}
&\qquad&
(1,1,0) & \mapsto &
\raisebox{-3pt}{
\setlength{\unitlength}{1.2pt}
\begin{picture}(60,10)
\put(0,10){\line(1,0){60}}
\put(0,0){\line(1,0){60}}
\put(0,0){\line(0,1){10}}
\put(10,0){\line(0,1){10}}
\put(20,0){\line(0,1){10}}
\put(30,0){\line(0,1){10}}
\put(40,0){\line(0,1){10}}
\put(50,0){\line(0,1){10}}
\put(60,0){\line(0,1){10}}
\put(0,0){\makebox(10,10){$\bigcirc$}}
\put(20,0){\makebox(10,10){$\bigcirc$}}
\end{picture}
}
\\
(1,0,1) & \mapsto &
\raisebox{-3pt}{
\setlength{\unitlength}{1.2pt}
\begin{picture}(60,10)
\put(0,10){\line(1,0){60}}
\put(0,0){\line(1,0){60}}
\put(0,0){\line(0,1){10}}
\put(10,0){\line(0,1){10}}
\put(20,0){\line(0,1){10}}
\put(30,0){\line(0,1){10}}
\put(40,0){\line(0,1){10}}
\put(50,0){\line(0,1){10}}
\put(60,0){\line(0,1){10}}
\put(0,0){\makebox(10,10){$\bigcirc$}}
\put(40,0){\makebox(10,10){$\bigcirc$}}
\end{picture}
}
&\qquad&
(0,2,0) & \mapsto &
\raisebox{-3pt}{
\setlength{\unitlength}{1.2pt}
\begin{picture}(60,10)
\put(0,10){\line(1,0){60}}
\put(0,0){\line(1,0){60}}
\put(0,0){\line(0,1){10}}
\put(10,0){\line(0,1){10}}
\put(20,0){\line(0,1){10}}
\put(30,0){\line(0,1){10}}
\put(40,0){\line(0,1){10}}
\put(50,0){\line(0,1){10}}
\put(60,0){\line(0,1){10}}
\put(20,0){\makebox(10,10){$\bigcirc$}}
\put(30,0){\makebox(10,10){$\bigcirc$}}
\end{picture}
}
\\
(0,1,1) & \mapsto &
\raisebox{-3pt}{
\setlength{\unitlength}{1.2pt}
\begin{picture}(60,10)
\put(0,10){\line(1,0){60}}
\put(0,0){\line(1,0){60}}
\put(0,0){\line(0,1){10}}
\put(10,0){\line(0,1){10}}
\put(20,0){\line(0,1){10}}
\put(30,0){\line(0,1){10}}
\put(40,0){\line(0,1){10}}
\put(50,0){\line(0,1){10}}
\put(60,0){\line(0,1){10}}
\put(20,0){\makebox(10,10){$\bigcirc$}}
\put(40,0){\makebox(10,10){$\bigcirc$}}
\end{picture}
}
&\qquad&
(0,0,2) & \mapsto &
\raisebox{-3pt}{
\setlength{\unitlength}{1.2pt}
\begin{picture}(60,10)
\put(0,10){\line(1,0){60}}
\put(0,0){\line(1,0){60}}
\put(0,0){\line(0,1){10}}
\put(10,0){\line(0,1){10}}
\put(20,0){\line(0,1){10}}
\put(30,0){\line(0,1){10}}
\put(40,0){\line(0,1){10}}
\put(50,0){\line(0,1){10}}
\put(60,0){\line(0,1){10}}
\put(40,0){\makebox(10,10){$\bigcirc$}}
\put(50,0){\makebox(10,10){$\bigcirc$}}
\end{picture}
}
\end{array}
$$

Let $A = (a_{ij})_{1 \leq i \leq M, \, 1 \leq j \leq N}$ be any 
$M \times N$ matrix, and let $I = \{ i_{1} < \cdots < i_{r} \} \subseteq [M]$ 
(resp. $J = \{ j_{1} < \cdots < j_{r} \} \subseteq [N]$) be 
a row (resp. column) index set. 
Let $A^{I}_{J} = A^{i_{1},\dots,i_{r}}_{j_{1},\dots,j_{r}}$ 
denote the matrix obtained from $A$ by choosing rows indexed by $I$ 
and columns indexed by $J$. 
If $r = M$ and $I = [M]$ (i.e., we choose all rows), 
then we may write $A_{J} = A_{j_{1},\dots,j_{r}}$ for $A^{[M]}_{J}$ 
when there is no fear of confusion.
Similarly we may use the notation $A^{I} = A^{i_{1},\dots,i_{r}}$ 
for $A^{I}_{[N]}$ when $r = N$ and $J = [N]$.

The purpose of this paper is to prove the following theorem 
and give an application to identities for classical group characters.

\begin{theorem}
\label{th:main}
Let $s$ and $n$ be positive integers and
$A = ( a_{ij} )_{1 \le i \le s+n-1, \, 1 \le j \le sn}$ be
an $(s+n-1) \times sn$ matrix.
We put
$$
\mathcal{R} = \binom{[s+n-1]}{n},
\quad
\mathcal{C} = \{ \iota_{s,n}(\mu) \,|\, \mu \in Z_{s,n} \},
\quad
\mathcal{C}^{0} = \{ \iota_{s,n}(\mu) \,|\, \mu \in Z^0_{s,s+n-1} \}.
$$
(Hereafter we assume $\iota=\iota_{s,n}$
unless $s$ and $n$ are explicitly specified.)
Then we have
\begin{equation}
\det \left( \det A^{I}_{J} \right)_{I \in \mathcal{R}, \, J \in \mathcal{C}}
=
\prod_{J \in \mathcal{C}^{0}} \det A_{J},
\label{eq:main}
\end{equation}
where the rows and columns of the matrix on the left-hand side are
arranged in increasing order with respect to $<$.
\end{theorem}

For example, if $n=1$, then
$$
\mathcal{R} = \mathcal{C} = \{ \{ 1 \}, \dots, \{ s \} \},
\quad
\mathcal{C}^{0} = \{ [s] \}
$$
and Equation (\ref{eq:main}) trivially holds.
If $s=1$, then
$$
\mathcal{R} = \mathcal{C} = \{ [n] \},
\quad
\mathcal{C}^{0} = \{ [n] \}
$$
and Equation (\ref{eq:main}) is also trivial.
If $s=2$ and $n=2$, then Equation (\ref{eq:main}) reads
$$
\det \begin{pmatrix}
\det A^{12}_{12} & \det A^{12}_{13} & \det A^{12}_{34} \\\noalign{\smallskip}
\det A^{13}_{12} & \det A^{13}_{13} & \det A^{13}_{34} \\\noalign{\smallskip}
\det A^{23}_{12} & \det A^{23}_{13} & \det A^{23}_{34}
\end{pmatrix}
= \det A^{123}_{123} \cdot \det A^{123}_{134}.
$$
If $s=2$ and $n=3$, then Equation (\ref{eq:main}) reads
\begin{equation*}
\det \begin{pmatrix}
\det A^{123}_{123} & \det A^{123}_{124} & \det A^{123}_{145} & \det A^{123}_{456} \\\noalign{\smallskip}
\det A^{124}_{123} & \det A^{124}_{124} & \det A^{124}_{145} & \det A^{124}_{456} \\\noalign{\smallskip}
\det A^{134}_{123} & \det A^{134}_{124} & \det A^{134}_{145} & \det A^{134}_{456} \\\noalign{\smallskip}
\det A^{234}_{123} & \det A^{234}_{124} & \det A^{234}_{145} & \det A^{234}_{456}
\end{pmatrix}
=
\det A^{1234}_{1234} \cdot \det A^{1234}_{1245} \cdot
\det A^{1234}_{1456}.
\end{equation*}
If $s=3$ and $n=2$, then Equation (\ref{eq:main}) reads
\begin{multline*}
\det \begin{pmatrix}
\det A^{12}_{12} & \det A^{12}_{13} & \det A^{12}_{15}
& \det A^{12}_{34} & \det A^{12}_{35} & \det A^{12}_{56} \\\noalign{\smallskip}
\det A^{13}_{12} & \det A^{13}_{13} & \det A^{13}_{15}
& \det A^{13}_{34} & \det A^{13}_{35} & \det A^{13}_{56} \\\noalign{\smallskip}
\det A^{14}_{12} & \det A^{14}_{13} & \det A^{14}_{15}
& \det A^{14}_{34} & \det A^{14}_{35} & \det A^{14}_{56} \\\noalign{\smallskip}
\det A^{23}_{12} & \det A^{23}_{13} & \det A^{23}_{15}
& \det A^{23}_{34} & \det A^{23}_{35} & \det A^{23}_{56} \\\noalign{\smallskip}
\det A^{24}_{12} & \det A^{24}_{13} & \det A^{24}_{15}
& \det A^{24}_{34} & \det A^{24}_{35} & \det A^{24}_{56} \\\noalign{\smallskip}
\det A^{34}_{12} & \det A^{34}_{13} & \det A^{34}_{15}
& \det A^{34}_{34} & \det A^{34}_{35} & \det A^{34}_{56}
\end{pmatrix}
\\
=
\det A^{1234}_{1235} \cdot \det A^{1234}_{1345} \cdot \det A^{1234}_{1356}.
\end{multline*}
We note that, as polynomials in indeterminates $a_{ij}$'s, 
the degree of the left hand side of (\ref{eq:main}) is equal to
$$
n \cdot \# Z_{s,n} = n \cdot \binom{s+n-1}{n}
$$
and the degree of the right hand side is equal to
$$
(s+n-1) \cdot \# Z^0_{s,s+n-1} = (s+n-1) \cdot \binom{s+n-2}{n-1},
$$
so that the degrees coincide with each other.

\section{
A Proof of the main theorem
}

In this section we give a proof of the main theorem, 
i.e., Eq.~(\ref{eq:main}).
Our proof splits into the following three steps:
\begin{enumerate}
\item
Finding the irreducible factors of the determinant on the 
left-hand side of (\ref{eq:main}) 
(Lemma~\ref{lem:factor1}),
\item
Proving Eq.~(\ref{eq:main}) up to constant 
(Lemma~\ref{lem:factor2}),
\item
Determining the constant (Lemma~\ref{lem:constant}).
\end{enumerate}

\subsection{
A partial order on the set $Z_{s,n}$ of compositions
}

Let $s$ and $n$ be positive integers, and fix a positive integer $k$ 
such that $1 \leq k \leq s$. 
Hereafter we call $k$ a \newterm{color}
and write $C = \{ 1, 2, \dots, s \}$ (the set of colors).
We introduce a partial order $\preceq_{k}$ on the set $Z_{s,n}$ 
of compositions as follows.
For $\lambda$ and $\mu$ in $Z_{s,n}$, we define $\lambda \preceq_{k} \mu$ 
if $\lambda_{i} \leq \mu_{i}$ for all $i \neq k$.
For example, if $\lambda = (2,0,1,3)$, $\mu = (2,1,2,1) \in Z_{4,6}$ 
and $k=4$, then we have $\lambda \preceq_{4} \mu$ 
since $\lambda_{i} \leq \mu_{i}$ for $i \neq 4$.
This gives a partially ordered set $(Z_{s,n},\preceq_{k})$, 
and we write $\lambda \prec_{k} \mu$ if $\lambda \preceq_{k} \mu$ 
and $\lambda \neq \mu$. 
Throughout this section, we will follow the general notation 
concerning posets used in the book \cite{Sta}. 
We see that $(Z_{s,n},\preceq_{k})$ is a graded poset 
with the rank function 
$$
\rho^{(k)}(\mu)
= n - \mu_{k}
= \sum_{i \neq k} \mu_{i}
\quad
\text{for}
\quad
\mu = (\mu_{1},\dots,\mu_{s}) \in Z_{s,n}.
$$
This poset $(Z_{s,n}, \preceq_{k})$ has the minimum element 
$\mu_{0}^{(k)} = (0, \dots, 0, n, 0, \dots, 0)$ in which $n$ 
is at the $k$th position.
Let $P_{i}^{(k)} = \left\{ \mu \in Z_{s,n} \,|\, \rho^{(k)}(\mu) 
= i \right\}$ be the subset of elements of rank $i$.
By definition, $Z_{s,n} = P_{0}^{(k)} \sqcup \cdots \sqcup 
P_{n}^{(k)}$ is a disjoint union,
where $P_{0}^{(k)} = \{ \mu_{0}^{(k)} \}$ and $P_{n}^{(k)} 
= \left\{ \mu=(\mu_1, \dots, \mu_s) \in Z_{s,n} \,|\, \mu_k=0
\right\}$.
Note that $P_{n}^{(k)}$ can be naturally identified with $Z_{s-1,n}$.
Thus we have $\# P_{0}^{(k)}=1$ and $\# P_{n}^{(k)}=\binom{s+n-2}{n}$.
For example, 
the Hasse diagram of $(Z_{3,2},\preceq_{1})$
is shown as Figure~\ref{fig:{Z}_{3,2}},
and we have 
$P_{0}^{(1)} = \{ (2,0,0) \}$, 
$P_{1}^{(1)} = \{ (1,1,0), (1,0,1) \}$ 
and 
$P_{2}^{(1)} = \{ (0,2,0), (0,1,1), (0,0,2) \}$.
%
%
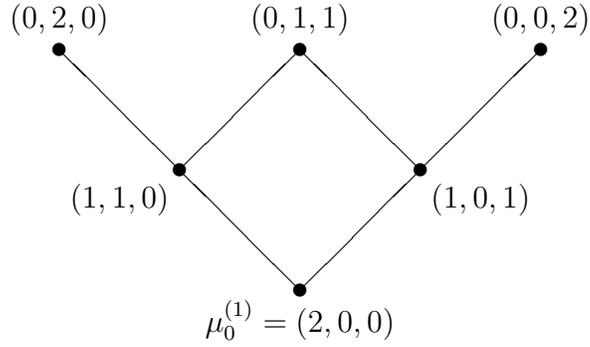
\begin{figure}[h] 
\begin{center}
\setlength{\unitlength}{0.8mm}
\begin{picture}(100,60)
\put( 50, 10){\line(-1, 1){20}}
\put( 50, 10){\line( 1, 1){20}}
\put( 30, 30){\line(-1, 1){20}}
\put( 30, 30){\line( 1, 1){20}}
\put( 70, 30){\line(-1, 1){20}}
\put( 70, 30){\line( 1, 1){20}}
\put( 10, 50){\circle*{2}}
\put( 50, 50){\circle*{2}}
\put( 90, 50){\circle*{2}}
\put( 30, 30){\circle*{2}}
\put( 70, 30){\circle*{2}}
\put( 50, 10){\circle*{2}}
\put( 10, 55){\makebox(0,0){$(0,2,0)$}}
\put( 50, 55){\makebox(0,0){$(0,1,1)$}}
\put( 90, 55){\makebox(0,0){$(0,0,2)$}}
\put( 20, 25){\makebox(0,0){$(1,1,0)$}}
\put( 80, 25){\makebox(0,0){$(1,0,1)$}}
\put( 50,  5){\makebox(0,0){$\mu^{(1)}_{0}=(2,0,0)$}}
\end{picture}
\caption{Hasse diagram of $\left({Z}_{3,2},\preceq_{1}\right)$}
\label{fig:{Z}_{3,2}}
\end{center}
\end{figure}
%
%

We put
$$
P^{(k)}
= \bigsqcup_{i=0}^{n-1}P_{i}^{(k)}
= Z_{s,n}\setminus P_{n}^{(k)},
$$
which is the set of elements $\mu = (\mu_{1}, \dots, \mu_{s})$ in $Z_{s,n}$
such that $\mu_{k}>0$.
Then we have $\# P^{(k)} = \binom{s+n-1}{n}-\binom{s+n-2}{n}
= \binom{s+n-2}{n-1}$.
We define a bijection $\tau^{(k)}: P^{(k)} \rightarrow Z^0_{s,s+n-1}$ 
as follows.
To each $\mu \in P^{(k)}$, we associate $\nu = (\nu_{1},\dots,\nu_{s}) = \tau^{(k)}(\mu)$ 
defined by
\begin{equation}
\nu_i = 
\begin{cases}
\mu_i   & \text{ if $i = k$,} \\
\mu_i+1 & \text{ if $i \neq k$.}
\end{cases}
\label{eq:tau}
\end{equation}
The inverse map of $\tau^{(k)}$ associates to $\nu \in Z^0_{s,s+n-1}$
a composition $\mu = (\mu_1, \dots, \mu_s)$ defined by
$$
\mu_{i} =
\begin{cases}
\nu_{i}   &\text{ if $i= k$,} \\
\nu_{i}-1 &\text{ if $i\neq k$.}
\end{cases}
$$
We also introduce another map $\phi^{(k)} : P^{(k)} \to \binom{[sn]}{s-1}$, 
which will play a crucial role in the following arguments.
To each $\mu \in P^{(k)}$, define $\phi^{(k)}(\mu) \in \binom{[sn]}{s-1}$ 
by
\begin{equation}
\phi^{(k)}(\mu)
= 
\bigsqcup_{\substack{1 \leq i \leq s \\ i \neq k}}
\{(i-1)n+\mu_{i}+1\}.
\label{eq:phi}
\end{equation}
In the above example where $s=3$, $n=2$ and $k=1$, we have
\begin{align*}
P^{(1)}
&= P_{0}^{(1)} \sqcup P_{1}^{(1)}
= \{ (2,0,0), (1,1,0), (1,0,1) \},
\\
Z^0_{3,4}
&= \{ (2,1,1), (1,2,1), (1,1,2) \}.
\end{align*}
The map $\tau^{(1)}$ defined in (\ref{eq:tau}) is given by
$$
\tau^{(1)}(2,0,0) = (2,1,1),
\quad
\tau^{(1)}(1,1,0) = (1,2,1),
\quad
\tau^{(1)}(1,0,1) = (1,1,2).
$$
And the map $\phi^{(1)} : P^{(1)} \rightarrow \binom{[6]}{2}$ 
defined in (\ref{eq:phi}) is given by
$$
\phi^{(1)}(2,0,0) = \{ 3, 5 \},
\quad
\phi^{(1)}(1,1,0) = \{ 4, 5 \},
\quad
\phi^{(1)}(1,0,1) = \{ 3, 6 \}.
$$
This map $\phi^{(1)}$ can be depicted as follows:
$$
\begin{array}{ccc}
& &
\raisebox{-3pt}{
\setlength{\unitlength}{1.2pt}
\begin{picture}(60,10)
\put(0,0){\makebox(10,10){$1$}}
\put(10,0){\makebox(10,10){$2$}}
\put(20,0){\makebox(10,10){$3$}}
\put(30,0){\makebox(10,10){$4$}}
\put(40,0){\makebox(10,10){$5$}}
\put(50,0){\makebox(10,10){$6$}}
\end{picture}
}
\\
(2,0,0) & \mapsto &
\raisebox{-3pt}{
\setlength{\unitlength}{1.2pt}
\begin{picture}(60,10)
\put(0,10){\line(1,0){60}}
\put(0,0){\line(1,0){60}}
\put(0,0){\line(0,1){10}}
\put(10,0){\line(0,1){10}}
\put(20,0){\line(0,1){10}}
\put(30,0){\line(0,1){10}}
\put(40,0){\line(0,1){10}}
\put(50,0){\line(0,1){10}}
\put(60,0){\line(0,1){10}}
\put(20,0){\makebox(10,10){$\bigcirc$}}
\put(40,0){\makebox(10,10){$\bigcirc$}}
\end{picture}
}
\\
(1,1,0) & \mapsto &
\raisebox{-3pt}{
\setlength{\unitlength}{1.2pt}
\begin{picture}(60,10)
\put(0,10){\line(1,0){60}}
\put(0,0){\line(1,0){60}}
\put(0,0){\line(0,1){10}}
\put(10,0){\line(0,1){10}}
\put(20,0){\line(0,1){10}}
\put(30,0){\line(0,1){10}}
\put(40,0){\line(0,1){10}}
\put(50,0){\line(0,1){10}}
\put(60,0){\line(0,1){10}}
\put(30,0){\makebox(10,10){$\bigcirc$}}
\put(40,0){\makebox(10,10){$\bigcirc$}}
\end{picture}
}
\\
(1,0,1) & \mapsto &
\raisebox{-3pt}{
\setlength{\unitlength}{1.2pt}
\begin{picture}(60,10)
\put(0,10){\line(1,0){60}}
\put(0,0){\line(1,0){60}}
\put(0,0){\line(0,1){10}}
\put(10,0){\line(0,1){10}}
\put(20,0){\line(0,1){10}}
\put(30,0){\line(0,1){10}}
\put(40,0){\line(0,1){10}}
\put(50,0){\line(0,1){10}}
\put(60,0){\line(0,1){10}}
\put(20,0){\makebox(10,10){$\bigcirc$}}
\put(50,0){\makebox(10,10){$\bigcirc$}}
\end{picture}
}
\end{array}
$$

Then one easily sees that the following proposition holds from the definition: 

\begin{prop}
\label{prop:bijection}
Fix a color $k \in C$.
For each $\mu\in P^{(k)}$, 
we have
$$
\iota(\mu) \cap \phi^{(k)}(\mu) = \emptyset,
\quad
\iota(\mu) \sqcup \phi^{(k)}(\mu)
= \iota \left( \tau^{(k)}(\mu) \right).
$$
\end{prop}

\subsection{
Vector notation
}

We write $|I|$ for $\sum_{i \in I} i$ if $I= \{ i_1, \dots, i_p \} 
\in \binom{[N]}{p}$.
If $I$ is a subset of $S$,
we denote the complement of $I$ in $S$ by $\overline{I}$ 
(or more explicitly by $S \setminus I$).

For two subset $I$ and $J$ of $[N]$, we define $\ep(I,J)
\in \{ 1, 0, -1 \}$ as follows.
If $I \cap J \neq \emptyset$, then we define $\ep(I,J) = 0$.
Otherwise, if $I = \{ i_1 < \cdots < i_p \}$ and $J = 
\{ j_1 < \cdots < j_q \}$, then $\ep(I,J)$ is the signature of 
a permutation which transforms the sequence $(i_1, \dots, i_p, 
j_1, \dots, j_q)$ in the increasing order.
For example, if $I = \{ 3, 4, 8 \}$ and $J = \{ 1, 6, 9 \}$,
then we have
$$
\ep(I,J)
=
\sgn
\begin{pmatrix}
1 & 3 & 4 & 6 & 8 & 9 \\
3 & 4 & 8 & 1 & 6 & 9
\end{pmatrix}
=
(-1)^4 = 1,
$$
where $4$ is the inversion number of the above permutation.

Let $s$ and $n$ be positive integers,
and let $A$ be an $(s+n-1)\times sn$ matrix.
For $J \in \binom{[sn]}{n}$ and $K \in \binom{[sn]}{s-1}$, 
we put
\begin{gather}
\mathcal{V}_{J}(A)
 =
\left( \det A^{I}_{J} \right)_{I \in \binom{[s+n-1]}{n}},
\label{eq:V}
\\
\overline{\mathcal{V}}_{K}(A)
 =
\left( (-1)^{|I|-n(n+1)/2} \det A^{\overline{I}}_{K} \right)_{I \in \binom{[s+n-1]}{n}},
\label{eq:barV}
\end{gather}
which are both $\binom{s+n-1}{n}$-dimensional column vectors, 
where the entries are arranged increasingly in the lexicographic 
order of indices.
For example,
if $s=3$, $n=2$, $J = \{ 1, 3 \} \in \binom{[6]}{2}$ and 
$K = \{ 4, 6 \} \in \binom{[6]}{2}$, then
\begin{gather*}
\mathcal{V}_{J}(A)
 = 
\trans\begin{pmatrix}
\det A^{12}_{13} & \det A^{13}_{13} & \det A^{14}_{13}
 & \det A^{23}_{13} & \det A^{24}_{13} & \det A^{34}_{13}
\end{pmatrix},
\\
\overline{\mathcal{V}}_{K}(A)
 = 
\trans\begin{pmatrix}
\det A^{34}_{46} & - \det A^{24}_{46} & \det A^{23}_{46}
 & \det A^{14}_{46} & - \det A^{13}_{46} & \det A^{12}_{46}
\end{pmatrix}.
\end{gather*}

For two vectors $v = \left( v_I \right)_{I \in \binom{[s+n-1]}{n}}$,
 $w = \left( w_I \right)_{I \in \binom{[s+n-1]}{n}}$, we write
$$
\langle v, w \rangle
 = \trans v w
 = \sum_{I \in \binom{[s+n-1]}{n}} v_I w_I.
$$
Using this inner-product notation, the Laplace expansion formula 
can be written as
\begin{equation}
\left\langle 
\mathcal{V}_{J}(A), \overline{\mathcal{V}}_{K}(A) 
\right\rangle
 =
\ep(J, K) \det A_{J \sqcup K},
\label{eq:Laplace}
\end{equation}
for $J \in \binom{[sn]}{n}$ and $K \in \binom{[sn]}{s-1}$.
If we take the above example $J = \{ 1, 3 \} \in \binom{[6]}{2}$ 
and $K = \{ 4, 6 \} \in \binom{[6]}{2}$, 
then we have
\begin{align*}
\left\langle 
\mathcal{V}_{J}(A), \overline{\mathcal{V}}_{K}(A) 
\right\rangle
&=
 \det A^{12}_{13}\det A^{34}_{46}
-\det A^{13}_{13}\det A^{24}_{46}
+\det A^{14}_{13}\det A^{23}_{46}
\\
&\quad
+\det A^{23}_{13}\det A^{14}_{46}
-\det A^{24}_{13}\det A^{13}_{46}
+\det A^{34}_{13}\det A^{12}_{46}
\\
&=
 \det A^{1234}_{1346}.
\end{align*}
The following proposition will be used to compute the determinant 
in this section.

\begin{prop}
\label{prop:inner-prod}
Let $s$ and $n$ be positive integers,
and let $A$ be an $(s+n-1)\times sn$ matrix.
Fix a color $k \in C$.
Let $\lambda\in Z_{s,n}$ and $\mu\in P^{(k)}$.
Then we have
$$
\left\langle 
\mathcal{V}_{\iota(\lambda)}(A), \overline{\mathcal{V}}_{\phi^{(k)}(\mu)}(A) 
\right\rangle
=0,
$$
unless $\lambda \preceq_{k} \mu$.
\end{prop}

\begin{demo}{Proof}
Assume $\lambda \not\preceq_{k} \mu$.
Then there exists $l \neq k$ such that $\lambda_l > \mu_l$.
Since this implies $(l-1)n+\mu_{l}+1 \in \iota(\lambda) \cap \phi^{(k)}(\mu)$, 
we obtain $\iota(\lambda) \cap \phi^{(k)}(\mu) \neq \emptyset$.
Thus we have
$\left\langle 
 \mathcal{V}_{\iota(\lambda)}(A),\overline{\mathcal{V}}_{\phi^{(k)}(\mu)}(A) 
\right\rangle=0$
from (\ref{eq:Laplace}).
\end{demo}

Let $s$ and $n$ be positive integers,
and let $A=(a_{ij})_{1 \leq i \leq s+n-1, \, 1 \leq j \leq sn}$ 
be an $(s+n-1) \times sn$ matrix.
We arrange the column vectors $\mathcal{V}_{\iota(\mu)}(A)$ 
($\mu\in Z_{s,n}$) into an $\binom{s+n-1}{n} \times \binom{s+n-1}{n}$ 
matrix $\mathcal{M}(A)$ :
\begin{equation}
\mathcal{M}(A)
 =
\left(
 \mathcal{V}_{\iota(\mu)}(A)
\right)_{\mu \in Z_{s,n}}
 =
\left(
 \det A^{I}_{\iota(\mu)}
\right)_{I \in \binom{[s+n-1]}{n},\, \mu \in Z_{s,n}},
\label{eq:M(A)}
\end{equation}
where the indices $\mu$ are arranged in the increasing order with respect to $<$.
For example,
if $s = 3$ and $n = 2$,
then
$$
\mathcal{M}(A)
 =
\begin{pmatrix}
\det A^{12}_{12} & \det A^{12}_{13} & \det A^{12}_{15}
 & \det A^{12}_{34} & \det A^{12}_{35} & \det A^{12}_{56} \\
\noalign{\smallskip}
\det A^{13}_{12} & \det A^{13}_{13} & \det A^{13}_{15}
 & \det A^{13}_{34} & \det A^{13}_{35} & \det A^{13}_{56} \\
\noalign{\smallskip}
\det A^{14}_{12} & \det A^{14}_{13} & \det A^{14}_{15}
 & \det A^{14}_{34} & \det A^{14}_{35} & \det A^{14}_{56} \\
\noalign{\smallskip}
\det A^{23}_{12} & \det A^{23}_{13} & \det A^{23}_{15}
 & \det A^{23}_{34} & \det A^{23}_{35} & \det A^{23}_{56} \\
\noalign{\smallskip}
\det A^{24}_{12} & \det A^{24}_{13} & \det A^{24}_{15}
 & \det A^{24}_{34} & \det A^{24}_{35} & \det A^{24}_{56} \\
\noalign{\smallskip}
\det A^{34}_{12} & \det A^{34}_{13} & \det A^{34}_{15}
 & \det A^{34}_{34} & \det A^{34}_{35} & \det A^{34}_{56}
\end{pmatrix}.
$$
Then our goal Eq.~(\ref{eq:main}) is to prove
\begin{equation}
\det \mathcal{M}(A)
=
\prod_{J \in \mathcal{C}^{0}}
\det A_{J},
\label{eq:main2}
\end{equation}
where $\mathcal{C}^{0} = \{ \iota(\mu) \,|\, \mu \in Z^0_{s,s+n-1} \}$.

Given an injection 
\begin{equation}
\Phi : Z_{s,n} \rightarrow \binom{[sn]}{s-1},
\label{eq:Phi}
\end{equation}
we consider an $\binom{s+n-1}{n} \times \binom{s+n-1}{n}$ matrix
defined by
\begin{equation}
\widehat{\mathcal{M}}(\Phi,A)
 =
\left(
 \overline{\mathcal{V}}_{\Phi(\mu)}(A)
\right)_{\mu \in Z_{s,n}}
 =
\left(
 (-1)^{|I|-n(n+1)/2} \det A^{\overline I}_{\Phi(\mu)}
\right)_{I \in \binom{[s+n-1]}{n}, \, \mu \in Z_{s,n}},
\label{eq:Phi-A}
\end{equation}
where the indices $\mu$ are arranged increasingly with respect to the order defined above.
Then, from \eqref{eq:Laplace}, we obtain
\begin{align}
\det \mathcal{M}(A) \cdot \det \widehat{\mathcal{M}}(\Phi,A)
 &=
\det \left(
 \trans\mathcal{M}(A) \widehat{\mathcal{M}}(\Phi,A)
\right)
\nonumber \\
 &=
\det \left(
 \left\langle
 \mathcal{V}_{\iota(\lambda)}(A), \overline{\mathcal{V}}_{\Phi(\mu)}(A) 
 \right\rangle
\right)_{\lambda, \, \mu \in Z_{s,n}}
\nonumber \\
 &=
\det \left(
 \ep( \iota(\lambda), \Phi(\mu) )
 \det A_{ \iota(\lambda) \sqcup \Phi(\mu) }
\right)_{\lambda, \, \mu \in Z_{s,n}}.
\label{eq:det-inprod}
\end{align}
In the next subsection 
we choose appropriate maps $\Phi$ which enable us to compute the last determinant 
in (\ref{eq:det-inprod}).
Since $\Rat[a_{ij} \,|\, 1 \leq i \leq s+n-1, 1 \leq j \leq sn]$ 
is a unique factorization domain, 
the irreducible factors of $\det \mathcal{M}(A)$ 
must appear in the last determinant.

\subsection{
Finding factors of the determinant
}

The goal of this subsection is to prove \eqref{eq:main2} up to constant.
In this subsection, we work in the polynomial ring 
$S = \Rat[a_{i,j} \,|\, 1 \leq i \leq s+n-1, 1 \leq j \leq sn]$ 
over the rational number field $\Rat$ in the indeterminates $a_{ij}$'s.

First we prove the following lemma which says 
each irreducible factor of $\det \mathcal{M}(A)$ is 
of the form $\det A_{\iota(\nu)}$, $\nu \in Z^0_{s,s+n-1}$.

\begin{lemma}
\label{lem:factor1}
Let $s$ and $n$ be positive integers.
Let $A$ be an $(s+n-1) \times sn$ matrix.
Then there exist non-negative integers $m_{\nu}$, $\nu \in Z^0_{s,s+n-1}$,
and a constant $c \in \Rat$ such that
\begin{equation}
\det \mathcal{M}(A)
 =
c
\prod_{\nu \in Z^0_{s,s+n-1}}
\left( \det A_{\iota(\nu)} \right)^{m_\nu}
\label{eq:lem1}
\end{equation}
holds.
\end{lemma}

\begin{demo}{Proof}
Let $\pi: Z_{s,n} \to C = \{ 1, 2, \dots, s \}$ be a map 
determined by the following condition:
\begin{quote}
For $\mu \in Z_{s,n}$, the color $k = \pi(\mu)$ is the least index
satisfying $\mu_{k} = \max \{ \mu_{l} \,|\, l = 1, \dots, s \}$.
\end{quote}
Note that, if $\pi(\mu) = k$, then $\mu_k > 0$ and $\mu \in P^{(k)}$.

Define a map $\Phi$ in \eqref{eq:Phi} by
$$
\Phi(\mu) = \phi^{(\pi(\mu))}(\mu)
\quad\text{for}\quad \mu \in Z_{s,n},
$$
where $\phi^{(k)}$ is as defined in \eqref{eq:phi}.
Then we claim that the determinant in \eqref{eq:det-inprod} equals
\begin{equation}
\det \left(
\left\langle
 \mathcal{V}_{\iota(\lambda)}(A), \overline{\mathcal{V}}_{\Phi(\mu)}(A) 
\right\rangle
\right)_{\lambda, \, \mu \in Z_{s,n}}
 =
\pm
\prod_{\mu \in Z_{s,n}} \det A_{\iota(\mu) \sqcup \Phi(\mu)}.
\label{eq:claim-det1}
\end{equation}
If we prove \eqref{eq:claim-det1}, 
then we obtain the desired result (\ref{eq:lem1}) 
since $\iota(\mu) \sqcup \phi^{(k)}(\mu) \in \iota(Z^0_{s,s+n-1})$ 
for any $k$ from Proposition~\ref{prop:bijection},
and $\det A_{J}$, $J \in \binom{[sn]}{s+n-1}$, are irreducible polynomials
in the unique factorization domain $S$.

To prove \eqref{eq:claim-det1}, 
set $m = m_{s,n} = n - \lceil n/s \rceil$,
where $\lceil x \rceil$ denotes the greatest integer which is 
less than or equal to $x$. 
Then $m$ satisfies $0 \leq m < n$.
We set
$$
Q_i =
\left\{
 \mu = (\mu_{1}, \dots, \mu_{s}) \in Z_{s,n} \,|\,
 \max \{ \mu_{l} \,|\, l = 1, \dots, s \} = n-i
\right\}
\quad(0 \leq i \leq m).
$$
Then we have 
$Z_{s,n} = Q_{0} \sqcup Q_{1} \sqcup \cdots \sqcup Q_{m}$.
For example, if $s = 3$ and $n = 2$, then we have
$m = m_{3,2} = 2 - \lceil 2/3 \rceil = 1$, 
$Q_{0} = \left\{ (2,0,0), (0,2,0), (0,0,2) \right\}$ and
$Q_{1} = \left\{ (1,1,0), (1,0,1), (0,1,1) \right\}$.
Then we claim that 
the determinant on the left-hand side of \eqref{eq:claim-det1}
 satisfies
\begin{multline}
\det \left(
 \left\langle
  \mathcal{V}_{\iota(\lambda)}(A),\overline{\mathcal{V}}_{\Phi(\mu)}(A) 
 \right\rangle
\right)_{\lambda, \mu \in Z_{s,n}}
\\
=
\pm
\prod_{\mu \in Q_{<r}}
 \det A_{\iota(\mu) \sqcup \Phi(\mu)}
\cdot
\det \left(
 \left\langle
 \mathcal{V}_{\iota(\lambda)}(A),\overline{\mathcal{V}}_{\Phi(\mu)}(A)
 \right\rangle
\right)_{\lambda, \, \mu \in Q_{\geq r}}
\label{eq:key1}
\end{multline}
for $0 \leq r \leq m+1$, where
$$
Q_{<r} = \bigsqcup_{i=0}^{r-1} Q_i,
\quad
Q_{\geq r} = \bigsqcup_{i=r}^{m} Q_i.
$$
We prove this identity (\ref{eq:key1}) by induction on $r$.
If $r=0$ then (\ref{eq:key1}) is trivial.
Assume (\ref{eq:key1}) is true for $r=t$.
If $\lambda \in Q_{\geq t}$ and $\mu \in Q_{t}$, 
then $\lambda_k \leq n-t$ and $\mu_{k} = n-t$ where $k = \pi(\mu)$.
Thus, if 
$
\left\langle
 \mathcal{V}_{\iota(\lambda)}(A),\overline{\mathcal{V}}_{\phi^{\left(\pi(\mu)\right)}(\mu)}(A) 
\right\rangle
\neq 0$, 
then 
$\lambda \preceq_{k} \mu$ from Proposition~\ref{prop:inner-prod},
which implies $\lambda = \mu$.
By expanding the determinant 
$
\det \left(
 \left\langle
 \mathcal{V}_{\iota(\lambda)}(A),\overline{\mathcal{V}}_{\Phi(\mu)}(A) 
 \right\rangle
\right)_{\lambda, \, \mu \in Q_{\geq t}}
$
along $\mu$th column for $\mu \in Q_{t}$, 
we obtain that this determinant equals
$$
\pm
\prod_{\mu \in Q_{t}}
 \det A_{\iota(\mu) \sqcup \Phi(\mu)}
\cdot
\det\left(
 \left\langle
  \mathcal{V}_{\iota(\lambda)}(A),\overline{\mathcal{V}}_{\Phi(\mu)}(A) 
 \right\rangle
\right)_{\lambda, \, \mu \in Q_{\geq t+1}},
$$
which implies (\ref{eq:key1}) is true when $r=t+1$.
Considering the case $r=m+1$ in (\ref{eq:key1}), 
we obtain \eqref{eq:claim-det1}.
This proves our lemma.
\end{demo}

\begin{example}
\label{ex:matrix1}
If $s = 3$ and $n = 2$, 
then the map $\pi : Z_{3,2} \rightarrow C = \{ 1, 2 \}$ 
is given by
\begin{gather*}
\pi(2,0,0) = 1, \quad \pi(1,1,0) = 1, \quad \pi(1,0,1) = 1,
\\
\pi(0,2,0) = 2, \quad \pi(0,1,1) = 2, \quad \pi(0,0,2) = 3.
\end{gather*}
Thus,
the map $\Phi : Z_{3,2} \to \binom{[6]}{2}$ ($\Phi(\mu) =
 \phi^{(\pi(\mu))}(\mu)$) is visualized as follows:
$$
\begin{array}{ccccccc}
 & &
\raisebox{-3pt}{
\setlength{\unitlength}{1.2pt}
\begin{picture}(60,10)
\put(0,0){\makebox(10,10){$1$}}
\put(10,0){\makebox(10,10){$2$}}
\put(20,0){\makebox(10,10){$3$}}
\put(30,0){\makebox(10,10){$4$}}
\put(40,0){\makebox(10,10){$5$}}
\put(50,0){\makebox(10,10){$6$}}
\end{picture}
}
&\qquad&
 & &
\raisebox{-3pt}{
\setlength{\unitlength}{1.2pt}
\begin{picture}(60,10)
\put(0,0){\makebox(10,10){$1$}}
\put(10,0){\makebox(10,10){$2$}}
\put(20,0){\makebox(10,10){$3$}}
\put(30,0){\makebox(10,10){$4$}}
\put(40,0){\makebox(10,10){$5$}}
\put(50,0){\makebox(10,10){$6$}}
\end{picture}
}
\\
(2,0,0) & \mapsto &
\raisebox{-3pt}{
\setlength{\unitlength}{1.2pt}
\begin{picture}(60,10)
\put(0,10){\line(1,0){60}}
\put(0,0){\line(1,0){60}}
\put(0,0){\line(0,1){10}}
\put(10,0){\line(0,1){10}}
\put(20,0){\line(0,1){10}}
\put(30,0){\line(0,1){10}}
\put(40,0){\line(0,1){10}}
\put(50,0){\line(0,1){10}}
\put(60,0){\line(0,1){10}}
\put(20,0){\makebox(10,10){$\bigcirc$}}
\put(40,0){\makebox(10,10){$\bigcirc$}}
\end{picture}
}
&\qquad&
(1,1,0) & \mapsto &
\raisebox{-3pt}{
\setlength{\unitlength}{1.2pt}
\begin{picture}(60,10)
\put(0,10){\line(1,0){60}}
\put(0,0){\line(1,0){60}}
\put(0,0){\line(0,1){10}}
\put(10,0){\line(0,1){10}}
\put(20,0){\line(0,1){10}}
\put(30,0){\line(0,1){10}}
\put(40,0){\line(0,1){10}}
\put(50,0){\line(0,1){10}}
\put(60,0){\line(0,1){10}}
\put(30,0){\makebox(10,10){$\bigcirc$}}
\put(40,0){\makebox(10,10){$\bigcirc$}}
\end{picture}
}
\\
(1,0,1) & \mapsto &
\raisebox{-3pt}{
\setlength{\unitlength}{1.2pt}
\begin{picture}(60,10)
\put(0,10){\line(1,0){60}}
\put(0,0){\line(1,0){60}}
\put(0,0){\line(0,1){10}}
\put(10,0){\line(0,1){10}}
\put(20,0){\line(0,1){10}}
\put(30,0){\line(0,1){10}}
\put(40,0){\line(0,1){10}}
\put(50,0){\line(0,1){10}}
\put(60,0){\line(0,1){10}}
\put(20,0){\makebox(10,10){$\bigcirc$}}
\put(50,0){\makebox(10,10){$\bigcirc$}}
\end{picture}
}
&\qquad&
(0,2,0) & \mapsto &
\raisebox{-3pt}{
\setlength{\unitlength}{1.2pt}
\begin{picture}(60,10)
\put(0,10){\line(1,0){60}}
\put(0,0){\line(1,0){60}}
\put(0,0){\line(0,1){10}}
\put(10,0){\line(0,1){10}}
\put(20,0){\line(0,1){10}}
\put(30,0){\line(0,1){10}}
\put(40,0){\line(0,1){10}}
\put(50,0){\line(0,1){10}}
\put(60,0){\line(0,1){10}}
\put(0,0){\makebox(10,10){$\bigcirc$}}
\put(40,0){\makebox(10,10){$\bigcirc$}}
\end{picture}
}
\\
(0,1,1) & \mapsto &
\raisebox{-3pt}{
\setlength{\unitlength}{1.2pt}
\begin{picture}(60,10)
\put(0,10){\line(1,0){60}}
\put(0,0){\line(1,0){60}}
\put(0,0){\line(0,1){10}}
\put(10,0){\line(0,1){10}}
\put(20,0){\line(0,1){10}}
\put(30,0){\line(0,1){10}}
\put(40,0){\line(0,1){10}}
\put(50,0){\line(0,1){10}}
\put(60,0){\line(0,1){10}}
\put(0,0){\makebox(10,10){$\bigcirc$}}
\put(50,0){\makebox(10,10){$\bigcirc$}}
\end{picture}
}
&\qquad&
(0,0,2) & \mapsto &
\raisebox{-3pt}{
\setlength{\unitlength}{1.2pt}
\begin{picture}(60,10)
\put(0,10){\line(1,0){60}}
\put(0,0){\line(1,0){60}}
\put(0,0){\line(0,1){10}}
\put(10,0){\line(0,1){10}}
\put(20,0){\line(0,1){10}}
\put(30,0){\line(0,1){10}}
\put(40,0){\line(0,1){10}}
\put(50,0){\line(0,1){10}}
\put(60,0){\line(0,1){10}}
\put(0,0){\makebox(10,10){$\bigcirc$}}
\put(20,0){\makebox(10,10){$\bigcirc$}}
\end{picture}
}
\end{array}
$$
Then the $6 \times 6$ matrix $\widehat{\mathcal{M}}(\Phi,A)$ 
in (\ref{eq:Phi-A}) is given by
$$
\widehat{\mathcal{M}}(\Phi,A)
 =
\begin{pmatrix}
\det A^{34}_{35} & \det A^{34}_{45} & \det A^{34}_{36}
 & \det A^{34}_{15} & \det A^{34}_{16} & \det A^{34}_{13} \\
\noalign{\smallskip}
- \det A^{24}_{35} & - \det A^{24}_{45} & - \det A^{24}_{36}
 & - \det A^{24}_{15} & - \det A^{24}_{16} & - \det A^{24}_{13} \\
\noalign{\smallskip}
\det A^{23}_{35} & \det A^{23}_{45} & \det A^{23}_{36}
 & \det A^{23}_{15} & \det A^{23}_{16} & \det A^{23}_{13} \\
\noalign{\smallskip}
\det A^{14}_{35} & \det A^{14}_{45} & \det A^{14}_{36}
 & \det A^{14}_{15} & \det A^{14}_{16} & \det A^{14}_{13} \\
\noalign{\smallskip}
- \det A^{13}_{35} & - \det A^{13}_{45} & - \det A^{13}_{36}
 & - \det A^{13}_{15} & - \det A^{13}_{16} & - \det A^{13}_{13} \\
\noalign{\smallskip}
\det A^{12}_{35} & \det A^{12}_{45} & \det A^{12}_{36}
 & \det A^{12}_{15} & \det A^{12}_{16} & \det A^{12}_{13}
\end{pmatrix}.
$$
Hence we obtain that
\begin{align*}
&
\trans\mathcal{M}(A) \widehat{\mathcal{M}}(\Phi,A)
\\
&\quad=
\left(
 \left\langle
  \mathcal{V}_{\iota(\lambda)}(A),\overline{\mathcal{V}}_{\Phi(\mu)}(A) 
 \right\rangle
\right)_{\lambda, \, \mu \in Z_{3,2}}
\\
&\quad=
\begin{pmatrix}
 \det A_{1235} & \det A_{1245} & \det A_{1236}
 &             0 &             0 &             0 \\
\noalign{\smallskip}
             0 & \det A_{1345} &             0
 &             0 &             0 &             0 \\
\noalign{\smallskip}
             0 &             0 &-\det A_{1356}
 &             0 &             0 &             0 \\
\noalign{\smallskip}
             0 &             0 &             0
 & \det A_{1345} & \det A_{1346} &             0 \\
\noalign{\smallskip}
             0 &             0 &             0
 &             0 & \det A_{1356} &             0 \\
\noalign{\smallskip}
             0 &             0 &             0
 &             0 &             0 & \det A_{1356} \\
\end{pmatrix},
\end{align*}
whose determinant is equal to
$$
- \det A_{1235} \, (\det A_{1345})^2 \, (\det A_{1356})^3.
$$
\end{example}

\begin{lemma}
\label{lem:factor2}
Let $s$ and $n$ be positive integers.
Let $A$ be an $(s+n-1) \times sn$ matrix.
Then there is  a constant $c' \in \Rat$ such that
\begin{equation}
\det \mathcal{M}(A)
 =
c'
\prod_{\nu \in Z^0_{s,s+n-1}}
\det A_{\iota(\nu)}.
\label{eq:lem2}
\end{equation}
\end{lemma}

\begin{demo}{Proof}
We may assume $n\geq2$ without loss of generality.
First, we claim that
the integers $m_{\nu}$, $\nu \in Z^0_{s,n}$, in (\ref{eq:lem1}) 
must be $0$ or $1$.
To prove this claim, 
we choose another map $\Phi$ in (\ref{eq:Phi}) 
and use the corresponding $\binom{s+n-1}{n} \times \binom{s+n-1}{n}$ matrix 
$\widehat{\mathcal{M}}(\Phi,A)$ defined in (\ref{eq:Phi-A}).
Fix a color $k_0 \in C$.
In the first step of our proof we construct 
a map $\Phi: Z_{s,n} \rightarrow \binom{[sn]}{s-1}$ 
whose restriction to $P^{(k_0)} = Z_{s,n} \setminus P^{(k_0)}_{n}$ 
agrees with $\phi^{(k_0)}$ given by (\ref{eq:phi}) 
and which satisfies the following three conditions:
\begin{enumerate}
\item[(i)]\label{cond:phi1}
for any $\lambda$, $\mu \in Z_{s,n}$, 
we have $\left\langle 
 \mathcal{V}_{\iota(\lambda)}(A),\overline{\mathcal{V}}_{\Phi(\mu)}(A) 
\right\rangle=0$
unless $\lambda \preceq_{k_0} \mu$,
\item[(ii)]\label{cond:phi2}
$\iota(\mu) \cap \Phi(\mu) = \emptyset$ for $\mu \in Z_{s,n}$,
\item[(iii)]\label{cond:phi3}
if $\mu \in P^{(k_0)}_{n}$, 
then $\iota(\mu) \sqcup \Phi(\mu) \not \in Z^0_{s,s+n-1}$.
\end{enumerate}
To construct such a map $\Phi$, we divide $P_{n}^{(k_0)}$ into two subsets.
Recall that $\mu_{0}^{(l)}$ is the composition whose $l$th entry is $n$ and others are $0$.
If we put $R^{(k_0)} = \{ \mu_{0}^{(l)} \,|\, l \neq k_0 \}$, 
then $R^{(k_0)}$ is an $(s-1)$-element subset of $P_{n}^{(k_0)}$.
We set the value $\Phi(\mu)$ for $\mu \in P_{n}^{(k_0)}$ as follows:
\begin{enumerate}
\item[(a)]
if $\mu \in P_{n}^{(k_0)} \setminus R^{(k_0)}$, 
we define $\Phi(\mu)=\phi^{(k_0)}(\mu)$ by (\ref{eq:phi}),
i.e.
\begin{equation*}
\Phi(\mu)
 =
 \bigsqcup_{\substack{1 \leq i \leq s \\ i \neq k_0}}
 \{ (i-1)n+\mu_{i}+1 \},
\end{equation*}
\item[(b)]
if $\mu_{0}^{(l)} \in R^{(k_0)}$ ($l \neq k_0$), 
we set
\begin{equation}
\Phi \left( \mu_{0}^{(l)} \right)
 =
\bigsqcup_{\substack{1 \leq i \leq s \\ i \neq k_0, \, l}}
 \{ (i-1)n + 1 \}
 \sqcup \{ (k_0-1)n+2 \}.
\label{eq:def_Phi}
\end{equation}
\end{enumerate}

To see that the condition (i) is fulfilled, 
it is enough to consider the case when $\mu \in P^{(k_0)}_{n}$.
If $\mu \in P_{n}^{(k_0)} \setminus R^{(k_0)}$,
then the proof is exactly the same argument as we did in the proof of Proposition~\ref{prop:inner-prod}.
If $\lambda \not\preceq_{k_0} \mu_{0}^{(l)}$ ($l \neq k_0$), 
then there exists $i \neq k_0, l$ such that $\lambda_i > 0$.
Thus we conclude from (\ref{eq:def_Phi}) that $\Phi( \mu_{0}^{(l)} ) \cap \iota(\lambda) 
\neq \emptyset$,
which immediately implies 
$\left\langle 
\mathcal{V}_{\iota(\lambda)}(A),\overline{\mathcal{V}}_{\Phi(\mu)}(A) 
\right\rangle
= 0$
by (\ref{eq:Laplace}).
This proves (i) for $\mu \in R^{(k_0)}$.
The condition (ii) is almost obvious from the definition.
The condition (iii) is also fulfilled 
since $(k_0-1)n+1 \not\in \Phi(\mu)$ and $(k_0-1)n+2 \in \Phi(\mu)$ 
for $\mu \in P^{(k_0)}_{n}$.

Now we claim that the determinant in \eqref{eq:det-inprod} equals
\begin{equation}
\det \left(
 \left\langle
  \mathcal{V}_{\iota(\lambda)}(A),\overline{\mathcal{V}}_{\Phi(\mu)}(A) 
 \right\rangle
\right)_{\lambda, \, \mu \in Z_{s,n}}
=
\pm
\prod_{\mu \in Z_{s,n}}
 \det A_{\iota(\mu) \sqcup \Phi(\mu)}.
\label{eq:claim-det2}
\end{equation}
It follows from Proposition~\ref{prop:bijection} that
\begin{equation*}
\prod_{\mu \in P_{<n}^{(k_0)}}
 \det A_{\iota(\mu) \sqcup \Phi(\mu)}
=
\prod_{\nu \in Z^0_{s,s+n-1}}
 \det A_{\iota(\nu)},
\end{equation*}
where
$$
P^{(k_0)}_{<n} = \bigsqcup_{i=0}^{n-1} P^{(k_0)}_i,
$$
and from (ii) and (iii) that 
$\iota(\mu) \sqcup \Phi(\mu) \not\in Z^0_{s,s+n-1}$
and $\det A_{\iota(\mu) \sqcup \Phi(\mu)} \neq 0$ 
for $\mu \in P_{n}^{(k_0)}$.
Thus, if we prove \eqref{eq:claim-det2}, 
then we can conclude that 
the integers $m_{\nu}$, $\nu \in Z^0_{s,s+n-1}$, in (\ref{eq:lem1}) 
must be $0$ or $1$ 
since $S = \Rat[a_{ij}]$ is a unique factorization domain.

To prove \eqref{eq:claim-det2}, 
it is enough to prove that 
\begin{multline}
\det \left(
 \left\langle
  \mathcal{V}_{\iota(\lambda)}(A),\overline{\mathcal{V}}_{\Phi(\mu)}(A) 
 \right\rangle
\right)_{\lambda, \, \mu \in Z_{s,n}}
\\
=
\pm
\prod_{\mu \in P_{<r}^{(k_0)}}
 \det A_{\iota(\mu) \sqcup \Phi(\mu)}
\cdot
\det\left(
 \left\langle
  \mathcal{V}_{\iota(\lambda)}(A),\overline{\mathcal{V}}_{\Phi(\mu)}(A) 
 \right\rangle
\right)_{\lambda, \, \mu \in P^{(k_0)}_{\geq r}}
\label{eq:key2}
\end{multline}
holds for $0 \leq r \leq n+1$, where
$$
P^{(k_0)}_{<r} = \bigsqcup_{i=0}^{r-1} P^{(k_0)}_i,
\quad
P^{(k_0)}_{\geq r} = \bigsqcup_{i=r}^{n} P^{(k_0)}_i.
$$
But, this can be proven by induction on $r$ 
exactly in the same way as we proved \eqref{eq:key1}
in the proof of Lemma~\ref{lem:factor1}.

Lastly we need to show that $m_{\nu}=1$ for all $\nu \in Z^0_{s,n}$.
In fact one easily sees this by comparing the degree of 
the both sides of (\ref{eq:lem1}).
The left-hand side $\det \mathcal{M}(A)$ is a polynomial of degree 
$n \binom{s+n-1}{n}$,
while the right-hand side is a polynomial of degree, 
at most, $(n+s-1) \binom{s+n-2}{n-1}$.
Thus all $m_{\nu}$ must be $1$.
This proves our lemma.
\end{demo}

\begin{example}
\label{ex:matrix2}
If $s = 3$, $n = 2$ and $k_0 = 1$,
then we have 
$$
P^{(1)} = \{ (2,0,0), (1,1,0), (1,0,1) \},
\quad
R^{(1)} = \{ (0,2,0), (0,0,2) \},
$$
and the above $\Phi$ takes values:
\begin{gather*}
\Phi(2,0,0) = \{ 3, 5 \}, \quad
\Phi(1,1,0) = \{ 4, 5 \}, \quad
\Phi(1,0,1) = \{ 3, 6 \},
\\
\Phi(0,2,0) = \{ 2, 5 \}, \quad
\Phi(0,1,1) = \{ 4, 6 \}, \quad
\Phi(0,0,2) = \{ 2, 3 \}.
\end{gather*}
$$
\begin{array}{ccccccc}
 & &
\raisebox{-3pt}{
\setlength{\unitlength}{1.2pt}
\begin{picture}(60,10)
\put(0,0){\makebox(10,10){$1$}}
\put(10,0){\makebox(10,10){$2$}}
\put(20,0){\makebox(10,10){$3$}}
\put(30,0){\makebox(10,10){$4$}}
\put(40,0){\makebox(10,10){$5$}}
\put(50,0){\makebox(10,10){$6$}}
\end{picture}
}
&\qquad&
 & &
\raisebox{-3pt}{
\setlength{\unitlength}{1.2pt}
\begin{picture}(60,10)
\put(0,0){\makebox(10,10){$1$}}
\put(10,0){\makebox(10,10){$2$}}
\put(20,0){\makebox(10,10){$3$}}
\put(30,0){\makebox(10,10){$4$}}
\put(40,0){\makebox(10,10){$5$}}
\put(50,0){\makebox(10,10){$6$}}
\end{picture}
}
\\
(2,0,0) & \mapsto &
\raisebox{-3pt}{
\setlength{\unitlength}{1.2pt}
\begin{picture}(60,10)
\put(0,10){\line(1,0){60}}
\put(0,0){\line(1,0){60}}
\put(0,0){\line(0,1){10}}
\put(10,0){\line(0,1){10}}
\put(20,0){\line(0,1){10}}
\put(30,0){\line(0,1){10}}
\put(40,0){\line(0,1){10}}
\put(50,0){\line(0,1){10}}
\put(60,0){\line(0,1){10}}
\put(20,0){\makebox(10,10){$\bigcirc$}}
\put(40,0){\makebox(10,10){$\bigcirc$}}
\end{picture}
}
&\qquad&
(1,1,0) & \mapsto &
\raisebox{-3pt}{
\setlength{\unitlength}{1.2pt}
\begin{picture}(60,10)
\put(0,10){\line(1,0){60}}
\put(0,0){\line(1,0){60}}
\put(0,0){\line(0,1){10}}
\put(10,0){\line(0,1){10}}
\put(20,0){\line(0,1){10}}
\put(30,0){\line(0,1){10}}
\put(40,0){\line(0,1){10}}
\put(50,0){\line(0,1){10}}
\put(60,0){\line(0,1){10}}
\put(30,0){\makebox(10,10){$\bigcirc$}}
\put(40,0){\makebox(10,10){$\bigcirc$}}
\end{picture}
}
\\
(1,0,1) & \mapsto &
\raisebox{-3pt}{
\setlength{\unitlength}{1.2pt}
\begin{picture}(60,10)
\put(0,10){\line(1,0){60}}
\put(0,0){\line(1,0){60}}
\put(0,0){\line(0,1){10}}
\put(10,0){\line(0,1){10}}
\put(20,0){\line(0,1){10}}
\put(30,0){\line(0,1){10}}
\put(40,0){\line(0,1){10}}
\put(50,0){\line(0,1){10}}
\put(60,0){\line(0,1){10}}
\put(20,0){\makebox(10,10){$\bigcirc$}}
\put(50,0){\makebox(10,10){$\bigcirc$}}
\end{picture}
}
&\qquad&
(0,2,0) & \mapsto &
\raisebox{-3pt}{
\setlength{\unitlength}{1.2pt}
\begin{picture}(60,10)
\put(0,10){\line(1,0){60}}
\put(0,0){\line(1,0){60}}
\put(0,0){\line(0,1){10}}
\put(10,0){\line(0,1){10}}
\put(20,0){\line(0,1){10}}
\put(30,0){\line(0,1){10}}
\put(40,0){\line(0,1){10}}
\put(50,0){\line(0,1){10}}
\put(60,0){\line(0,1){10}}
\put(10,0){\makebox(10,10){$\bigcirc$}}
\put(40,0){\makebox(10,10){$\bigcirc$}}
\end{picture}
}
\\
(0,1,1) & \mapsto &
\raisebox{-3pt}{
\setlength{\unitlength}{1.2pt}
\begin{picture}(60,10)
\put(0,10){\line(1,0){60}}
\put(0,0){\line(1,0){60}}
\put(0,0){\line(0,1){10}}
\put(10,0){\line(0,1){10}}
\put(20,0){\line(0,1){10}}
\put(30,0){\line(0,1){10}}
\put(40,0){\line(0,1){10}}
\put(50,0){\line(0,1){10}}
\put(60,0){\line(0,1){10}}
\put(30,0){\makebox(10,10){$\bigcirc$}}
\put(50,0){\makebox(10,10){$\bigcirc$}}
\end{picture}
}
&\qquad&
(0,0,2) & \mapsto &
\raisebox{-3pt}{
\setlength{\unitlength}{1.2pt}
\begin{picture}(60,10)
\put(0,10){\line(1,0){60}}
\put(0,0){\line(1,0){60}}
\put(0,0){\line(0,1){10}}
\put(10,0){\line(0,1){10}}
\put(20,0){\line(0,1){10}}
\put(30,0){\line(0,1){10}}
\put(40,0){\line(0,1){10}}
\put(50,0){\line(0,1){10}}
\put(60,0){\line(0,1){10}}
\put(10,0){\makebox(10,10){$\bigcirc$}}
\put(20,0){\makebox(10,10){$\bigcirc$}}
\end{picture}
}
\end{array}
$$
Thus the $6 \times 6$ matrix $\widehat{\mathcal{M}}(\Phi,A)$ 
in (\ref{eq:Phi-A}) equals
$$
\begin{pmatrix}
 \det A^{34}_{35} & \det A^{34}_{45} & \det A^{34}_{36}
 & \det A^{34}_{25} & \det A^{34}_{46} &  \det A^{34}_{23} \\
\noalign{\smallskip}
- \det A^{24}_{35} & - \det A^{24}_{45} & - \det A^{24}_{36}
 & - \det A^{24}_{25} & - \det A^{24}_{46} & - \det A^{24}_{23} \\
\noalign{\smallskip}
 \det A^{23}_{35} & \det A^{23}_{45} & \det A^{23}_{36}
 & \det A^{23}_{25} & \det A^{23}_{46} &  \det A^{23}_{23} \\
\noalign{\smallskip}
 \det A^{14}_{35} & \det A^{14}_{45} & \det A^{14}_{36}
 & \det A^{14}_{25} & \det A^{14}_{46} &  \det A^{14}_{23} \\
\noalign{\smallskip}
- \det A^{13}_{35} & - \det A^{13}_{45} & - \det A^{13}_{36}
 & - \det A^{13}_{25} & - \det A^{13}_{46} & - \det A^{13}_{23} \\
\noalign{\smallskip}
 \det A^{12}_{35} & \det A^{12}_{45} & \det A^{12}_{36}
 & \det A^{12}_{25} & \det A^{12}_{46} &  \det A^{12}_{23}
\end{pmatrix}.
$$
Hence we obtain that
\begin{align*}
&
\trans\mathcal{M}(A) \widehat{\mathcal{M}}(\Phi,A)
\\
&\quad=
\det \left(
 \left\langle
  \mathcal{V}_{\iota(\lambda)}(A),\overline{\mathcal{V}}_{\Phi(\mu)}(A) 
 \right\rangle
\right)_{\lambda, \, \mu \in Z_{s,n}}
\\
&\quad=
\begin{pmatrix}
\det A_{1235} & \det A_{1245} & \det A_{1236}
 &             0 & \det A_{1246} &             0 \\
\noalign{\smallskip}
            0 & \det A_{1345} &             0
 &-\det A_{1235} & \det A_{1346} &             0 \\
\noalign{\smallskip}
            0 &             0 &-\det A_{1356}
 &             0 &-\det A_{1456} & \det A_{1235} \\
\noalign{\smallskip}
            0 &             0 &             0
 & \det A_{2345} &             0 &             0 \\
\noalign{\smallskip}
            0 &             0 &             0
 &             0 &-\det A_{3456} &             0 \\
\noalign{\smallskip}
            0 &             0 &             0
 &             0 &             0 & \det A_{2356} 
\end{pmatrix},
\end{align*}
whose determinant is
$$
\det A_{2345} \, \det A_{3456} \, \det A_{2356} \, \prod_{\nu \in Z^0_{3,4}} \det A_{\iota(\nu)}.
$$
\end{example}

\subsection{
Determining the constant
}

Now we are in position to prove that $c'=1$ in Eq.~(\ref{eq:lem2}).
For this purpose we specialize the entries of
$A = ( a_{ij} )_{1 \le i \le s+n-1, \ 1 \le j \le sn}$ 
as
\begin{equation}
a_{i,j} = x_{j}^{s+n-i}
\quad(1 \le i \le s+n-1,\, 1 \le j \le sn).
\label{eq:special0}
\end{equation}
Here we introduce the \newterm{lexicographical monomial order} 
(or \newterm{lex order}) on the set of monomials 
in the variables $x_{1} \gg \cdots \gg x_{sn}$,
i.e.
$
x_{1}^{i_{1}} \cdots x_{sn}^{i_{sn}}
\gg
x_{1}^{j_{1}} \cdots x_{sn}^{j_{sn}}
$
if $i_{1} = j_{1}, \dots, i_{k-1} = j_{k-1}$ and $i_{k} > j_{k}$ 
for some $1 \leq k \leq sn$.
If $f$ is a polynomial of $x_{1}, \dots, x_{sn}$, 
then we define the \newterm{leading coefficient} of $f$ to be 
the coefficient of the greatest monomial in the lex order.
We also define the \newterm{leading term} of $f$ 
to be the leading monomial including its coefficient.
For example, 
the leading term of 
$f = 4 \, x_{1}x_{2}^2x_{3} + 4 \, x_{3}^2 - 5 \,x_{1}^3 x_{2} + 7 \, x_{1}^2 x_{3}^2$ 
is $- 5 \, x_{1}^3 x_{2}$ and its leading coefficient $-5$.
If $f$ and $g$ are two polynomials with the same leading term,
we shall write $f \sim g$.
In the above example, 
we can write $f \sim - 5\,x_{1}^3x_{2}$.
If $I = \{ i_{1} < \cdots < i_{n} \} \in \mathcal{R} 
= \binom{[s+n-1]}{n}$ 
and $J = \{ j_{1} < \cdots < j_{n} \} \in \mathcal{C} 
= \{ \iota(\mu) \,|\, \mu \in Z_{s,n} \}$, 
then the leading term of the determinant
\begin{equation*}
\det A^{I}_{J}
 =
\det \left( x_{j_{q}}^{s+n-i_{p}} \right)_{1 \leq p, \, q \leq n}
 =
\det
\begin{pmatrix}
 x_{j_{1}}^{s+n-i_{1}} & \hdots & x_{j_{n}}^{s+n-i_{1}} \\
 \vdots                & \ddots & \vdots \\
 x_{j_{1}}^{s+n-i_{n}} & \hdots & x_{j_{n}}^{s+n-i_{n}}
\end{pmatrix}
\end{equation*}
is obviously the product of its diagonal entries 
which equals $x_{j_{1}}^{s+n-i_{1}} \cdots x_{j_{n}}^{s+n-i_{n}}$.
Hereafter we shall write $x_{J}^{s+n-I}$ for 
$x_{j_{1}}^{s+n-i_{1}} \cdots x_{j_{n}}^{s+n-i_{n}}$.
Because of 
$\det A^{I}_{J} \sim x_{J}^{s+n-I}$, 
we conclude that
$$
\det \left( \det A^{I}_{J} \right)_{I \in \mathcal{R}, \, J \in \mathcal{C}}
\sim
\det \left( x_{J}^{s+n-I} \right)_{I \in \mathcal{R}, \, J \in \mathcal{C}}.
$$
Thus, if we prove the following lemma, 
then it shows that the leading coefficient of 
$\det \left( \det  A^{I}_{J} \right)_{I \in \mathcal{R}, \, J \in \mathcal{C}}$ 
equals $1$ under the specialization (\ref{eq:special0}).

\begin{lemma}
\label{lem:constant}
The leading term of 
$\det\left(x_{J}^{s+n-I}\right)_{I\in\mathcal{R},\,J\in\mathcal{C}}$
is given by the product of its diagonal entries. 
More explicitly, we have
\begin{equation}
\det \left(
 x_{J}^{s+n-I}
\right)_{I \in \mathcal{R}, \, J \in \mathcal{C}}
\sim
\prod_{k=1}^{s} \prod_{j=1}^{n}
 x_{(k-1)n+j}^{(s+1-k)\binom{s+n-j}{s}},
\label{eq:lt}
\end{equation}
and the leading coefficient is $1$.
\end{lemma}

\begin{demo}{Proof}
We proceed by double induction on $s$ and $n$.
In this proof we write $\mathcal{R}_{s,n}$ for $\mathcal{R}$ 
and $\mathcal{C}_{s,n}$ for $\mathcal{C}$.
If $s=1$ then the determinant equals
$
x_{1}^{n} x_{2}^{n-1} \cdots x_{n}
$,
and if $n=1$ then we have
$$
\det \left( x_{J}^{s+n-I} \right)_{I \in \mathcal{R}, \, J \in \mathcal{C}}
 =
\det \left( x_j^{s+n-i} \right)_{1 \leq i, \, j \leq s}
 \sim
x_{1}^{s} x_{2}^{s-1} \cdots x_{s}.
$$
This implies that our claim is valid for $s=1$ or $n=1$.

Assume $s > 1$ and $n > 1$.
We put 
\begin{alignat*}{2}
\mathcal{R}_{s,n}^{(0)}
 &=
\left\{
 \{ 1 \} \cup I \,\Bigl|\,
 I \in \binom{[2,s+n-1]}{n-1}
\right\},
&\qquad
\mathcal{R}_{s,n}^{(1)}
 &=
\binom{[2,s+n-1]}{n},
\\
\mathcal{C}_{s,n}^{(0)}
 &=
\iota(Z_{s,n}^{(0)}),
&\qquad
\mathcal{C}_{s,n}^{(1)}
 &=
\iota(Z_{s,n}^{(1)}),
\end{alignat*}
where
\begin{align*}
Z_{s,n}^{(0)}
 &=
\left\{ 
\mu = (\mu_{1}, \dots, \mu_{s}) \in Z_{s,n} \,|\, \mu_{1} > 0 \right\},
\\
Z_{s,n}^{(1)}
 &=
\left\{
\mu = (\mu_{1}, \dots, \mu_{s}) \in Z_{s,n} \,|\, \mu_{1} = 0 \right\}.
\end{align*}
Then we have 
$\mathcal{R}_{s,n} = \mathcal{R}_{s,n}^{(0)} \sqcup \mathcal{R}_{s,n}^{(1)}$,
$\mathcal{C}_{s,n} = \mathcal{C}_{s,n}^{(0)} \sqcup \mathcal{C}_{s,n}^{(1)}$,
$\# \mathcal{R}_{s,n}^{(0)} = \# \mathcal{C}_{s,n}^{(0)} = \binom{s+n-2}{n-1}$
and $\# \mathcal{R}_{s,n}^{(1)} = \# \mathcal{C}_{s,n}^{(1)} = \binom{s+n-2}{n}$.
Further, under the lexicographic order we defined in Section~1, 
the subset $\mathcal{R}_{s,n}^{(0)}$ consists of the first $\binom{s+n-2}{n-1}$ 
elements of $\mathcal{R}_{s,n}$, 
and $\mathcal{C}_{s,n}^{(0)}$ consists of the first $\binom{s+n-2}{n-1}$ 
elements of $\mathcal{C}_{s,n}$.
Expanding the determinant along the first $\binom{s+n-2}{n-1}$ columns, 
we obtain
\begin{align*}
\det\left(
 x_{J}^{s+n-I}
\right)_{I \in \mathcal{R}_{s,n},\,J \in \mathcal{C}_{s,n}}
=
\sum_{\mathcal{I}}
 \ep_{\mathcal{I}}
 \det \left(
  x_{J}^{s+n-I}
 \right)_{I \in \mathcal{I}, \, J \in \mathcal{C}_{s,n}^{(0)}}
 \det \left(
  x_{J}^{s+n-I}
 \right)_{I \in \overline{\mathcal{I}}, \, J \in \mathcal{C}_{s,n}^{(1)}},
\end{align*}
where $\mathcal{I}$ runs over all $\binom{s+n-2}{n-1}$-element 
subsets of $\mathcal{R}_{s,n}$, and $\ep_{\mathcal{I}} = \pm 1$.
The symbol $\overline{\mathcal{I}}$ stands for the complement set 
of $\mathcal{I}$ in $\mathcal{R}_{s,n}$.
Since the variable $x_{1}$ does not appear in 
$\det \left( x_{J}^{s+n-I} \right)_{I \in \overline{\mathcal{I}}, 
\, J \in \mathcal{C}_{s,n}^{(1)}}$, 
the exponent of $x_1$ attains its maximum when we take 
$\mathcal{I} = \mathcal{R}_{s,n}^{(0)}$, 
in which case we have $\ep_{\mathcal{R}_{s,n}^{(0)}}=1$.
Thus we obtain
$$
\det \left(
 x_{J}^{s+n-I}
\right)_{I \in \mathcal{R}_{s,n}, \, J \in \mathcal{C}_{s,n}}
\sim
\det \left(
 x_{J}^{s+n-I}
\right)_{I \in \mathcal{R}_{s,n}^{(0)}, \, J \in \mathcal{C}_{s,n}^{(0)}}
\det \left(
 x_{J}^{s+n-I}
\right)_{I \in \mathcal{R}_{s,n}^{(1)}, \, J \in \mathcal{C}_{s,n}^{(1)}}.
$$
Since we can naturally identify $\mathcal{R}_{s,n}^{(1)}$ and 
$\mathcal{C}_{s,n}^{(1)}$ with $\mathcal{R}_{s-1,n}$ and
$\mathcal{C}_{s-1,n}$ respectively, 
we obtain
\begin{equation}
\det \left(
 x_{J}^{s+n-I}
\right)_{I \in \mathcal{R}_{s,n}^{(1)}, \, J \in \mathcal{C}_{s,n}^{(1)}}
\sim
\det \left(
 x_{J}^{s+n-1-I}
\right)_{I \in \mathcal{R}_{s-1,n}, \, J \in \mathcal{C}_{s,n}^{(1)}}
\sim
\prod_{k=2}^{s} \prod_{j=1}^{n}
 x_{(k-1)n+j}^{(s+1-k)\binom{s+n-1-j}{s-1}}
\label{eq:lt1}
\end{equation}
by the induction hypothesis.
On the other hand, since each entry of
$
\det \left(
 x_{J}^{s+n-I}
\right)_{I \in \mathcal{R}_{s,n}^{(0)}, \, J \in \mathcal{C}_{s,n}^{(0)}}
$
has the factor $x_{1}^{s+n-1}$, 
we can sweep out the factor from each row 
and the remaining determinant has the form of type $(s,n-1)$ 
with the $s(n-1)$ variables $x_{2}, \dots, x_{n}$, $x_{n+2}, \dots,x_{2n}$, 
$\dots$, $x_{(s-1)n+2},\dots,x_{sn-1}$.
Thus we obtain
\begin{equation}
\det \left(
 x_{J}^{s+n-I}
\right)_{I \in \mathcal{R}_{s,n}^{(0)}, \, J \in \mathcal{C}_{s,n}^{(0)}}
\sim
\left(
 x_{1}^{s+n-1}
\right)^{\binom{s+n-2}{n-1}}
\times
\prod_{j=2}^{n}
 x_{j}^{s\binom{s+n-j}{s}}
\,
\prod_{k=2}^{s} \prod_{j=1}^{n-1}
 x_{(k-1)n+j}^{(s+1-k)\binom{s+n-1-j}{s}}.
\label{eq:lt2}
\end{equation}
It is an easy computation to check Eq.~(\ref{eq:lt}) by using 
(\ref{eq:lt1}) and (\ref{eq:lt2}).
This proves our lemma.
\end{demo}

Now we can finish our proof of Theorem~\ref{th:main}.

\begin{demo}{Proof of Theorem~\ref{th:main}}
Now we consider the right-hand side of Eq.~(\ref{eq:lem2})
under the specialization (\ref{eq:special0}).
For $J = \{ j_{1} < \cdots < j_{s+n-1} \} \in \mathcal{C}^{0}$, 
the leading term of $\det A_{J}$ is obviously the product 
of the diagonal entries, which shows
$$
\det A_{J}
\sim 
x^{s+n-[s+n-1]}_{J}
=
x_{j_{1}}^{s+n-1} x_{j_{2}}^{s+n-2} \cdots x_{j_{s+n-1}}.
$$
Thus it is obvious that 
the leading coefficient of $\prod_{\nu \in Z^0_{s,s+n-1}} \det A_{\iota(\nu)}$ 
equals $1$,
which implies $c'=1$ in (\ref{eq:lem2}).
In fact, by a direct computation, one can also obtain 
$$
\prod_{\nu \in Z^0_{s,s+n-1}} \det A_{\iota(\nu)}
 \sim
\prod_{k=1}^{s} \prod_{j=1}^{n}
 x_{(k-1)n+j}^{(s+1-k)\binom{s+n-j}{s}}.
$$
We complete the proof of the main formula (\ref{eq:main}).
\end{demo}

\section{
An application to determinants for classical group characters
}

In this section, we apply Theorem~\ref{th:main} to determinants 
involving classical group characters and Macdonald polynomials.
We almost follow the notation in \cite{IOW}.

A partition is a sequence $\lambda = (\lambda_1, \lambda_2, \dots)$ 
of non-negative integers in weakly decreasing order and 
containing only finitely many non-zero terms.
The non-zero $\lambda_i$ are call the \defterm{parts} of $\lambda$.
The number of parts is the \defterm{length} of $\lambda$, 
denoted by $\ell(\lambda)$,
and the sum of parts is the \defterm{weight} of $\lambda$, 
denoted by $|\lambda|$.
For two partitions $\lambda=(\lambda_{1}, \dots, \lambda_{m})$ 
and $\mu=(\mu_{1}, \dots, \mu_{m})$ with at most $m$ parts, 
we say $\lambda$ precedes $\mu$ in the \newterm{reverse lexicographic ordering} 
if the first non-vanishing difference $\lambda_{i}-\mu_{i}$ is positive.

If $k$ is a non-negative integer,
let $\delta (k) = (k,k-1,k-2, \dots, 1, 0, 0, \dots)$ denote
the staircase partition, 
and if $k$ is a positive half-integer,
let $\delta(k) = (k,k-1, \dots, 3/2, 1/2)$ denote
the decreasing sequence of length $(2k+1)/2$.
When $\vectx = (x_1, \dots, x_n)$ is an $n$-tuple of variables,
and $\alpha = (\alpha_1, \dots, \alpha_n)$ 
is a sequence of integers or a sequence of half-integers, 
we use the notation
$$
V(\alpha ;\vectx)
= \big( x_i^{\alpha_j} \big)_{1 \leq i, j \leq n},
\quad
W^{\pm}(\alpha ; \vectx)
= \big( x_i^{\alpha_j} \pm x_i^{-\alpha_j} \big)_{1 \leq i, j \leq n}.
$$
We now recall Weyl's character formula.
If $\lambda = (\lambda_1, \dots, \lambda_n)$ is a partition of length $\leq n$, 
we put
\begin{align}
\GL_n(\lambda;\vectx)
&=
\frac{ \det V(\lambda + \delta(n-1) ; \vectx) }
     { \det V(\delta(n-1) ; \vectx) },
\label{eq:char-gl}
\\
\Symp_{2n}(\lambda;\vectx)
&=
\frac{ \det W^-(\lambda + \delta(n) ; \vectx) }
     { \det W^-(\delta(n) ; \vectx) },
\label{eq:char-sp}
\\
\Orth_{2n+1}(\lambda;\vectx)
&=
\frac{ \det W^-(\lambda + \delta(n-1/2) ; \vectx) }
     { \det W^-(\delta(n-1/2) ; \vectx) },
\label{eq:char-o1}
\\
\Orth_{2n}(\lambda;\vectx)
&=
\begin{cases}
\dfrac{ \det W^+(\lambda + \delta(n-1) ; \vectx) }
      { (1/2) \det W^+(\delta(n-1) ; \vectx) }
&\text{if $\lambda_n \neq 0$,} \\
\dfrac{ \det W^+(\lambda + \delta(n-1) ; \vectx) }
      { \det W^+(\delta(n-1) ; \vectx) }
&\text{if $\lambda_n = 0$.} \\
\end{cases}
\label{eq:char-o2}
\end{align}
Then $\GL_n(\lambda;\vectx)$ (resp. $\Symp_{2n}(\lambda;\vectx)$) gives
the irreducible character of the general linear group $\GL_n$ 
(resp. the symplectic group $\Symp_{2n}$) with highest weight $\lambda$.
$\GL_n(\lambda;\vectx)$ is usually called the \newterm{Schur function}.
Similarly, $\Orth_N(\lambda;\vectx)$ is the irreducible character of 
the orthogonal group $\Orth_N$ with ``highest weight'' $\lambda$.
Weyl's denominator formula reads as follows.

\begin{prop}
\begin{gather*}
\det V(\delta(n-1) ; \vectx)
=
\prod_{1 \leq i < j \leq n} (x_{i}-x_{j}),
\\
\det W^-(\delta(n-1/2) ; \vectx)
=
(-1)^{n}
\prod_{i=1}^{n}
\frac{1-x_{i}}{x_{i}^{1/2}}
\prod_{1\leq i < j \leq n}
\frac{(x_{j}-x_{i})(1-x_{i}x_{j})}{x_{i}x_{j}},
\\
\det W^-(\delta(n) ; \vectx)
=
(-1)^{n}
\prod_{i=1}^{n}
\frac{1-x_{i}^2}{x_{i}}
\prod_{1 \leq i < j \leq n}
\frac{(x_{j}-x_{i})(1-x_{i}x_{j})}{x_{i}x_{j}},
\\
\det W^+(\delta(n-1) ; \vectx)
=
2
\prod_{1 \leq i < j \leq n}
\frac{(x_{j}-x_{i})(1-x_{i}x_{j})}{x_{i}x_{j}}.
\end{gather*}
\end{prop}

In what follows, we consider $sn$ variables $x^{(1)}_1, \dots, x^{(1)}_n,
x^{(2)}_1, \dots, x^{(2)}_n, \dots, x^{(s)}_1, \dots, x^{(s)}_n$, 
and we write
$$
X_\mu
= ( x^{(1)}_1, \dots, x^{(1)}_{\mu_1},
x^{(2)}_1, \dots, x^{(2)}_{\mu_2}, \dots,
x^{(s)}_1, \dots, x^{(s)}_{\mu_s}).
$$
for $\mu \in Z_{s,n}$,
We apply our main formula \eqref{eq:main} and obtain the following theorem.

\begin{theorem}
\label{th:Schur}
Let $s$ and $n$ be positive integers.
Then we have
\begin{equation}
\det \Big(
\GL_n(\lambda;X_\mu)
\Big)_{\lambda \subset ((s-1)^n), \, \mu \in  Z_{s,n}}
=
\prod_{1 \le k < l \le s} \prod_{i,j=1}^n
\left( x^{(k)}_i - x^{(l)}_j \right)^{\binom{s+n-i-j-1}{s-2}},
\label{eq:Cauchy-Sylvester-type-A}
\end{equation}
and
\begin{align}
&
\det \Big(
\Orth_{2n+1}(\lambda;X_\mu)
\Big)_{\lambda \subset ((s-1)^n), \, \mu \in  Z_{s,n}}
=
\det \Big(
\Symp_{2n}(\lambda;X_\mu)
\Big)_{\lambda \subset ((s-1)^n), \, \mu \in  Z_{s,n}}
\nonumber
\\
&
=
\det \Big(
\Orth_{2n}(\lambda;X_\mu)
\Big)_{\lambda \subset ((s-1)^n), \, \mu \in  Z_{s,n}}
\nonumber
\\
&
=
\prod_{1 \le k < l \le s} \prod_{i,j=1}^{n}
\left\{\frac
{  (x^{(l)}_j-x^{(k)}_i)(1-x^{(k)}_i  x^{(l)}_j)}
{ x^{(k)}_i x^{(l)}_j }
\right\}^{\binom{s+n-i-j-1}{s-2}},
\label{eq:Cauchy-Sylvester-type-BCD}
\end{align}
where the rows are indexed by partitions $\lambda
= (\lambda_1, \dots, \lambda_n)$ such that $\lambda_1 \le s-1$
and arranged increasingly in the reverse lexicographic ordering.
\end{theorem}

\begin{remark}
If we substitute 
$$
x^{(k)}_j = t^{j-1} a_k \quad(1 \le k \le s, \, 1 \le j \le n)
$$
into (\ref{eq:Cauchy-Sylvester-type-A}) and (\ref{eq:Cauchy-Sylvester-type-BCD}),
we obtain the (symplectic) Schur function identities in \cite[Corollary~3.1, Proposition~7.1]{AI2}.
\end{remark}

\begin{demo}{Proof of Theorem~\ref{th:Schur}}
We prove one case, say the symplectic case of 
(\ref{eq:Cauchy-Sylvester-type-BCD}), 
because the other cases can be proven similarly.
For simplicity, set
$$
C(u) = -\frac{1-u^2}{u},
\qquad
D(u,v) = \frac{ (v-u) (1-u v) }{ uv }.
$$
We apply the formula (\ref{eq:main}) to the matrix $A = (a_{ij})$ with entries
given by
\begin{multline*}
a_{i,(k-1)n+j}
=
\left( x^{(k)}_j \right)^{s+n-i}
-
\left( x^{(k)}_j \right)^{-s-n+i}
\\
(1 \le i \le s+n-1, \, 1 \le k \le s, \, 1 \le j \le n).
\end{multline*}
Then the $(I,J)$ entry of the matrix on the left-hand side of (\ref{eq:main})
is equal to
$$
\Delta_\mu \cdot \Symp_{2n}(\lambda;X_\mu)
$$
where $I$ and $J$ are given by
$$
I = \{ s-\lambda_1, s+1-\lambda_{2}, \dots, s+n-1-\lambda_{n} \},
\quad
J = \iota(\mu),
$$
and
$$
\Delta_\mu
=
\prod_{k=1}^{s} \prod_{i=1}^{\mu_{k}}C(x^{(k)}_i)\,
\prod_{k=1}^s \prod_{1 \le i < j \le \mu_k}
D( x^{(k)}_i , x^{(k)}_j )
\prod_{1 \le k < l \le s} \prod_{i=1}^{\mu_k} \prod_{j=1}^{\mu_l}
D( x^{(k)}_i , x^{(l)}_j).
$$
Hence we see that the left hand side of (\ref{eq:main}) becomes
$$
\prod_{\mu \in  Z_{s,n}} \Delta_\mu
\cdot
\det \Big(
\Symp_{2n}(\lambda; X_\mu )
\Big)_{\lambda \subset ((s-1)^n), \, \mu \in  Z_{s,n}}.
$$
The exponents of 
$C(x^{(k)}_i)$, $D(x^{(k)}_i ,x^{(k)}_j)$ ($i<j$) and $D(x^{(k)}_i , x^{(l)}_j)$ 
are given by 
\begin{gather*}
\# \{ \mu \in Z_{s,n} \,|\, \mu_k \ge i \}
= \# Z_{s,n-i}
= \binom{s+n-i-1}{n-i},
\\
\# \{ \mu \in  Z_{s,n} \,|\, \mu_k \ge j \}
= \#  Z_{s,n-j}
= \binom{s+n-j-1}{n-j},
\\
\# \{ \mu \in  Z_{s,n} \,|\, \mu_k \ge i, \ \mu_l \ge j \}
= \#  Z_{s,n-i-j}
= \binom{s+n-i-j-1}{n-i-j},
\end{gather*}
respectively.
On the other hand, the right-hand side of (\ref{eq:main}) becomes
$$
\prod_{\mu \in Z^0_{s,s+n-1}} \Delta_\mu,
$$
and the exponents of 
$C(x^{(k)}_i)$, $D(x^{(k)}_i ,x^{(k)}_j)$ ($i<j$) and $D(x^{(k)}_i , x^{(l)}_j)$ 
are given by
\begin{gather*}
\# \{ \mu \in Z^0_{s,s+n-1} \,|\, \mu_k \ge i \}
= \# Z^0_{s,s+n-1-(i-1)}
= \# Z^0_{s,s+n-i}
= \binom{s+n-i-1}{n-i},
\\
\# \{ \mu \in Z^0_{s,s+n-1} \,|\, \mu_k \ge j \}
= \# Z^0_{s,s+n-1-(j-1)}
= \# Z^0_{s,s+n-j}
= \binom{s+n-j-1}{n-j},
\\
\# \{ \mu \in Z^0_{s,s+n-1} \,|\, \mu_k \ge i, \ \mu_l \ge j \}
= \# Z^0_{s,s+n-1-(i-1)-(j-1)}
= \binom{s+n-i-j}{n-i-j+1},
\end{gather*}
respectively.
Therefore we see that
$
\det \Big(
\Symp_{2n}(\lambda; X_\mu )
\Big)_{\lambda \subset ((s-1)^n), \, \mu \in  Z_{s,n}}
$
equals
\begin{align*}
&
\prod_{k=1}^s \prod_{i=1}^{n}
C(x^{(k)}_i)
^{\binom{s+n-j-1}{n-j} - \binom{s+n-j-1}{n-j}}
\prod_{k=1}^s \prod_{1 \le i < j \le n}
D(x^{(k)}_i , x^{(k)}_j)^{\binom{s+n-j-1}{n-j} - \binom{s+n-j-1}{n-j}}
\\
&\quad\times
\prod_{1 \le k < l \le s} \prod_{i,j=1}^n
D(x^{(k)}_i , x^{(l)}_j)^{\binom{s+n-i-j}{n-i-j+1} - \binom{s+n-i-j-1}{n-i-j}}
\\
&=
\prod_{1 \le k < l \le s} \prod_{i,j=1}^n
D(x^{(k)}_i , x^{(l)}_j )^{\binom{s+n-i-j-1}{s-2}}.
\end{align*}
This proves the symplectic case.

One can prove the odd orthogonal case similarly 
by merely taking 
\begin{multline*}
a_{i,(k-1)n+j}
=
\left( x^{(k)}_j \right)^{s+n-i-1/2}
-
\left( x^{(k)}_j \right)^{-s-n+i+1/2}
\\
(1 \le i \le s+n-1, \, 1 \le k \le s, \, 1 \le j \le n),
\end{multline*}
and $C(u) = - u^{-1/2} (1-u)$ in the above arguments.

One should be more careful to prove the even orthogonal case 
because of the multiple $2$.
If one takes 
\begin{multline*}
a_{i,(k-1)n+j}
= 
\left( x^{(k)}_j \right)^{s+n-i-1}
+
\left( x^{(k)}_j \right)^{-s-n+i+1}
\\
(1 \le i \le s+n-1, \, 1 \le k \le s, \, 1 \le j \le n),
\end{multline*}
then the $(I,J)$ entry of the matrix on the left-hand side of (\ref{eq:main})
is 
$$
\Delta_\mu \cdot \Orth_{2n}(\lambda;X_\mu),
$$
where 
\begin{multline*}
\Delta_\mu
\\
=
\begin{cases}
\prod_{k=1}^s \prod_{1 \le i < j \le \mu_k}
D( x^{(k)}_i , x^{(k)}_j )
\prod_{1 \le k < l \le s} \prod_{i=1}^{\mu_k} \prod_{j=1}^{\mu_l}
D( x^{(k)}_i , x^{(l)}_j)
&\text{ if $\lambda_{n}\neq0$,}\\
2\prod_{k=1}^s \prod_{1 \le i < j \le \mu_k}
D( x^{(k)}_i , x^{(k)}_j )
\prod_{1 \le k < l \le s} \prod_{i=1}^{\mu_k} \prod_{j=1}^{\mu_l}
D( x^{(k)}_i , x^{(l)}_j)
&\text{ if $\lambda_{n}=0$,}
\end{cases}
\end{multline*}
which implies the left-hand side have the multiple $2^{\binom{s+n-2}{n-1}}$.
On the other hand, each determinant on the right-hand side has the factor $2$.
Hence $2^{\binom{s+n-2}{n-1}}$ also appears on the right-hand side.

Considering the principal term of asymptotic behavior of 
the both sides of (\ref{eq:Cauchy-Sylvester-type-BCD}) at $x_j^{(k)} \to +\infty$
($1 \le k \le s$, $1 \le j \le n$) 
(or just taking $a_{i,(k-1)n+j} = \left( x^{(k)}_j \right)^{s+n-i-1}$), 
we immediately obtain (\ref{eq:Cauchy-Sylvester-type-A}).
\end{demo}

Next we give a similar determinant identities for Macdonald polynomials.
Hereafter 
we also use a symmetric function $f$ in countably many variables 
$x=(x_1, x_2, \dots)$.
(See \cite{M} for detailed exposition on symmetric functions 
in countably many variables and Macdonald functions.)
The symmetric function $p_r = \sum_{i \geq 1} x_i^r$ is called 
the $r$th \newterm{power-sum symmetric function}.
The function
$$
m_{\lambda}(x) = \sum x^{\alpha}
$$
summed over all distinct permutations $\alpha$ of $\lambda$ is called 
the \newterm{monomial symmetric function}.
The \newterm{dominance order} is defined by
$$
\lambda \preceq \mu
\quad
\Longleftrightarrow
\quad
\lambda_{1}+\cdots+\lambda_{i} \leq \mu_{1}+\cdots+\mu_{i}
\text{ for all $i \geq 1$.}
$$

Define an inner product $\langle u, v \rangle_{q,t}$ 
by the values on the power-sum symmetric functions by
$$
\langle p_{\lambda}, p_{\mu} \rangle_{q,t}
=
\delta_{\lambda,\mu}
\,
z_{\lambda}
\,
\prod_{i=1}^{\ell(\lambda)}
\frac{1-q^{\lambda_{i}}}{1-t^{\lambda_{i}}},
$$
where $z_{\lambda} = \prod_{i \geq 1} i^{m_{i}} m_{i}!$
with $m_{i} = m_{i}(\lambda)$ the number of parts of $\lambda$ equal to $i$.
The Macdonald symmetric functions $P_{\lambda}(x;q,t)$ are defined 
by the property that they are upper triangularly related 
to the monomial symmetric functions and 
orthogonal with respect to the inner product. 
That is, they satisfy the following two conditions:
\begin{enumerate}
\item
$P_{\lambda}(x;q,t)
= m_{\lambda}(x) + \sum_{\mu \prec \lambda} c_{\lambda,\mu}(q,t) m_{\mu}(x)$, 
where the sum is over $\mu$ smaller than $\lambda$ in the dominance order.
\item
$\langle P_{\lambda}, P_{\mu} \rangle_{q,t}=0$ if $\lambda \neq \mu$.
\end{enumerate}
The Macdonald symmetric functions $Q_{\lambda}(x;q,t)$ are defined by
$$
Q_{\lambda}(x;q,t) = b_{\lambda}(q,t) P_{\lambda}(x;q,t)
$$
where
$$
b_{\lambda}(q,t)
=
\prod_{c \in\lambda} \frac{1-q^{a(c)}t^{l(c)+1}}{1-q^{a(c)+1}t^{l(c)}}
$$
with $a(c) = \lambda_{i}-j$ the arm-length and 
$l(c) = \lambda'_{j}-i$ the leg-length at $c = (i,j) \in \lambda$.
From Equation~(\ref{eq:Cauchy-Sylvester-type-A})
we immediately have the following:

\begin{corollary}
\label{cor;Macdonald}
Let $s$ and $n$ be positive integers.
Then we have
\begin{align}
\det \left(
P_\lambda( X_\mu;q,t )
\right)_{\lambda \subset ((s-1)^n), \, \mu \in  Z_{s,n}}
&=
\prod_{1 \le k < l \le s} \prod_{i,j=1}^n
\left( x^{(k)}_i - x^{(l)}_j \right)^{\binom{s+n-i-j-1}{s-2}},
\label{eq:MacdP}
\\
\det \left(
Q_\lambda( X_\mu;q,t )
\right)_{\lambda \subset ((s-1)^n), \, \mu \in  Z_{s,n}}
&=
\frac{(t^n;q)_{s-1}(t;t)_{n-1}}{(q;q)_{s-1}(tq^{s-1};t)_{n-1}}
\nonumber\\
&\quad\times
\prod_{1 \le k < l \le s} \prod_{i,j=1}^n
\left( x^{(k)}_i - x^{(l)}_j \right)^{\binom{s+n-i-j-1}{s-2}},
\label{eq:MacdQ}
\end{align}
where the rows are indexed by partitions $\lambda
= (\lambda_1, \dots, \lambda_n)$ such that $\lambda_1 \le s-1$
and arranged increasingly in the reverse lexicographic ordering.
\end{corollary}

\begin{demo}{Proof}
The first identity (\ref{eq:MacdP}) comes from the fact that
$$
P_{\lambda}(x;q,t)
=
\GL_n(\lambda;x) + \sum_{\mu \prec \lambda} a_{\lambda,\mu}(q,t) \GL_n(\mu;x)
$$
where the sum is over $\mu$ smaller than $\lambda$ in the dominance order.
This relation follows easily from the above definition and
$$
\GL_n(\lambda;x)
=
\sum_{\mu} K_{\lambda,\mu} m_{\mu}(x),
$$
where the Kostka number $K_{\lambda,\mu}$ satisfies $K_{\lambda,\mu}=0$ 
unless $\mu \preceq \lambda$ and $K_{\lambda,\lambda}=1$.
If $\nu \in Z_{s,n}$, 
then we have $\GL_n(\mu;X_\nu)=0$ unless $\ell(\mu) \leq n$.
Hence, for $\lambda \subset ((s-1)^n)$, and $\nu \in Z_{s,n}$,
we obtain
$$
P_{\lambda}(X_{\nu};q,t)
=
\GL_n(\lambda;X_\nu)
+
\sum_{\substack{\mu \subset ((s-1)^n) \\ \mu \prec \lambda}}
a_{\lambda,\mu}(q,t) \GL_n(\mu;X_\nu),
$$
which immediately implies (\ref{eq:MacdP}).
We also obtain the second identity (\ref{eq:MacdQ}) from the first one (\ref{eq:MacdP}) 
using
$$
\prod_{\lambda \subset ((s-1)^n)}b_{\lambda}(q,t)
=
\frac{(t^n;q)_{s-1}(t;t)_{n-1}}{(q;q)_{s-1}(tq^{s-1};t)_{n-1}}.
$$
This completes our proof.
\end{demo}

\end{document}